\documentstyle[11pt,epsf]{article}

\def\grr{\varrho}\def\fP{{\bf P}}\def\gD{\Delta}
\def\ga{\alpha}\def\fC{{\bf C}}\def\gS{\Sigma}\def\fZ{{\bf Z}}
\def\inn{\subset}\def\gk{\kappa}\def\cL{{\cal L}}
\def\bgm{\hbox{\boldmath$\mu$}}
\def\boldmu{{\bgm}} 
 \def\fR{{\bf R}}
\def\lra{\longrightarrow}\def\cF{{\cal F}}
\def\gs{\sigma}

\def\II{\rm II} \def\III{\rm III} \def\IV{\rm IV} \def\I{\rm I}
\newtheorem{theorem}{Theorem}[section]

\newtheorem{lemma}[theorem]{Lemma}

\begin{document}

\title{K3-fibered Calabi-Yau threefolds II, \\ singular fibers}
\author{Bruce Hunt } \maketitle

\bigskip
\begin{center}
  {\Large \bf Introduction}
\end{center}

In part I of this paper we have introduced a twist map from a product of
weighted hypersurfaces onto another weighted hypersurface. This map
realized a quotient formation. The image was shown to have the structure of
fibration. As it is known how to resolve singularities of these
hypersurfaces with the methods of toric geometry, this gives a method for
explicitly determining the singular fibers in that fibration. The twist map
is defined for a pair of weighted hypersurfaces $V_1\inn \fP_{(w_0,\ldots,
  w_n)}$ and $V_2 \inn \fP_{(v_0,\ldots, v_m)}$, and maps the product onto
a hypersurface in the weighted projective space $\fP_{(v_0 w_1, \ldots, v_0
  w_n, w_0 v_1, \ldots, w_0 v_m )}$. We are interested in the case that the
image is Calabi-Yau, i.e., has vanishing first Chern class, in which case
this Calabi-Yau has a constant-modulus fibration, either elliptic or K3.

For an introduction to this see part I. We will first recall the twist map
and some further notations, then begin by describing the Kodaira fibers of
types I$_0^*$, II, III, IV, II$^*$, III$^*$ and IV$^*$, from the same point
of view: as resolutions of images of the twist map. Here the reader will
see that this method is convenient and effective. Then we proceed to
consider the degererate fibers in families of K3 surfaces. Recall that
there is a vast liturature on stable degenerations of K3 surfaces. Looking
at elliptic curves, the stable degenerate ones are exactly the fibers of
type I$_n$, corresponding to a resolution of a surface $A_n$-singularity.
The singular fibers of types I$_0^*$, II, III, IV, II$^*$, III$^*$ and
IV$^*$ in Kodaira's classification, are all {\it non-stable}, and it is an
analog of these which we consider in this paper. 
We remark that in all cases we consider, the Calabi-Yau threefold has both
a K3- as well as an elliptic fibration. The elliptic fibrations are more
thoroughly studied, and we could have indeed considered these. However we
prefer the picture of K3-fibration for the following reason: all birational
transformations (i.e., non-uniqueness of models) take place in the fibers,
whereas for elliptic fibrations, to get certain models one must modify the
{\it base} of the fibration. Furthermore, this research may be considered a
first detailed look at non-stable K3-degenerations.

\section{Weighted projective spaces}
We will be working with weighted projective spaces, which are certain
(singular) quotients of usual projective space. Alternatively, they may be
described as quotients of $\fC^{n+1}$ by a $\fC^*$-action. We assume the
weights $(w_0,\ldots, w_n)$ are given, let $\boldmu_{w_i}$ denote the group
of $w_i$th roots of unity, and consider the action of $\boldmu:=
\boldmu_{w_1}\times \cdots \times \boldmu_{w_n}$ on $\fP^n$ as follows. Let
$g=(g_0,\ldots, g_n) \in \boldmu$, and consider for $(z_0:\ldots: z_n)$
homogenous coordinates on $\fP^n$ the action
\[ (g,(z_0:\ldots: z_n)) \mapsto (g_0z_0:\ldots: g_nz_n).\]
Alternatively, consider the action of $\fC^*$ on $\fC^{n+1}$ given by
\[ (t,(z_0,\ldots, z_n)) \mapsto (t^{w_0}z_0,\ldots, t^{w_n}z_n).\]
In both cases, the resulting quotient is the weighted projective space,
which we will denote by $\fP_{(w_0,\ldots, w_n)}$. General references for
weighted projective spaces are \cite{dolg} and \cite{Y}. A {\it weighted
  hypersurface} is the zero locus of a weighted homogenous polynomial $p$.
We will assume the weights are {\it normalized} in the sense that no $n$ of
the $n+1$ weights have a common divisor $>1$. Both for the weighted
projective spaces as well as for the weighted hypersurfaces this assumption
is no restriction (cf.\ \cite{dolg} 1.3.1 and \cite{Y}, pp.\ 185-186).  We
will write such isomorphisms in the sequel without further comment, for
example $\fP_{(2,3,6)}\cong \fP_{(2,1,2)}\cong \fP_{(1,1,1)}=\fP^2$, where
the first equality is because the last two weights are divisible by 3, the
second while the first and last are divisible by 2.

We will use the notation $\fP_{(w_0,\ldots, w_n)}[d]$ to denote either a
certain weighted hypersurface of degree $d$, or to denote the whole family
of such (the context will make the usage clear). In the particular case
that the weighted polynomial $p$ is of Fermat type, then there is a useful
fact, corresponding to the above normalizations. For example, in
$\fP_{(2,3,6)}$ consider the weighted hypersurface $x_0^6+x_1^4 +x_2^2=0$.
Then the isomorphism $\fP_{(2,3,6)}\cong \fP_{(2,1,2)}$ above is given by
the introduction of a new variable $(x_0')=x_0^3,$ which is in spite of
appearances a one to one coordinate transformation (becuase of admissible
rescalings), and the Fermat polynomial becomes $(x_0')^2+x_1^4+x_2^2=0$.
Again, the isomorphism $\fP_{(2,1,2)}\cong \fP_{(1,1,1)}$ is given by
setting $(x_1')=x_1^2$, and the Fermat polynomial becomes
$(x_0')^2+(x_1')^2+x_2^2=0$, which is a quadric in the projective plane. We
{\it denote} this process by the symbolic expressions
\[ \fP_{(2,3,6)}[12]\cong \fP_{(2,1,2)}[4] \cong \fP_{(1,1,1)}[2].\]

It is well-known how to resolve the weighted projective space
$\fP_{(w_0,\ldots, w_n)}$. For this, one takes the following vectors in
$\fR^n$,
\begin{equation}\label{vectors}
v_0 = {1\over w_0} \left(\begin{array}{c} -1 \\ -1 \\ \vdots \\ -1 
  \end{array}\right),\quad v_1={1\over w_1}\left(\begin{array}{c} 1 \\ 0
    \\ \vdots \\ 0  
  \end{array}\right), \quad \ldots, \quad v_n={1\over w_n} 
  \left(\begin{array}{c} 0 \\ \vdots \\ 0 \\ 1 
  \end{array}\right),
\end{equation}
and considers the lattice $\cL=\fZ v_0+\fZ v_1+\cdots + \fZ v_n$ in
$\fR^n$. The $n+1$ vectors $v_0,\ldots, v_n$ give a cone decomposition of
$\fR^n$, and this decomposition is {\it refined} until the resulting
decomposition satisfies: each cone has volume 1, where the volume is
normalized in such a way that the standard simplex in $\fR^n$ has volume
$\prod w_i$. This is equivalent to: any set of $n$ vectors spanning one of the
cones of the decomposition form a $\fZ$-basis of the lattice $\cL$.

Now suppose we are given a weighted hypersurface of degree $d$ in 
$\fP_{(w_0,\ldots, w_n)}$, such that $d=\sum w_i$. Then, as is well-known,
this is a sufficient condition for the variety to be Calabi-Yau, i.e., the
dualizing sheaf is trivial. Supposing moreover that the hypersurface is
quasi-smooth, then in dimensions 2 and 3, by work of Roan and Yau
(\cite{GRY}, section 3), 
there is a resolution of singularities such that the smooth
variety is still Calabi-Yau. In this case, the resolution (described in
\cite{RY}) is easier than of the ambient projective spaces themselves. The
reason is that it effectively reduces to a question of cones in one
dimension less. In particular, in the case of Calabi-Yau threefolds, the
resolution is described in terms of a simplicial decomposition of a
triangle (which is the {\it face} of one of the cones mentioned above). Let
$X \inn \fP_{\bf w}$ denote the singular weighted hypersurface. Then under
the assumption that $X$ is quasi-smooth, the singularities are all quotient
singularities by abelian groups. Locally they can be written as quotients
of $\fC^3$ by the following transformations
\begin{eqnarray*} \psi: \fC^3 & \lra & \fC^3 \\ 
(z_1,z_2,z_3) & \mapsto & \left(\exp({a\over d})z_1,\exp({b\over d})z_2,
\exp({c\over d})z_3\right),
\end{eqnarray*}
which describes an action of the group of $d$th roots of unity 
$\boldmu_d$ on $\fC^3$. Let $e^i,\ i=1,2,3$ denote the standard unit
vectors in $\fR^3$, and let 
\[ v=\left( \begin{array}{c} a/d \\ b/d \\ c/d
  \end{array} \right). \]
Let $\cL$ denote the lattice in $\fR^3$ spanned by the $e^i$ and
$v$. Finally, let $\gs$ denote the cone
\[ \gs=\left\{ \sum _{i=1} ^3 x_i e^i \in \fR^3 | x_i \geq 0,
  i=1,2,3\right\}. \]
Then this affine cone determines a toric variety, which is a neighborhood
of the singularity in question. Next one uses the fact that, since $X$ is
Calabi-Yau, $a+b+c$ is divisible by $d$, and this in turn implies that the
integral vectors we require to decompose the cone to get a smooth cone
decomposition all lie in a hypersurface. This is given by 
\[ \cF = \gs \cap \left\{ \sum x_i =1 \right\}.\]
This is a triangle, and we need to determine the number of vertices, edges
and two-simplices of the simplicial decomposition of $\cF$
to determine the number of resolution divisors, the
number of intersection curves and intersection points. 
First of all, since the area of $\cF$ is equal to $d$, there must be a
total of $s=d$ simplices in the decomposition. 
\begin{lemma}\label{formula} Define integers $d_i$ as follows.
\[ d_1=gcd(a,d),\quad d_2=gcd(b,d), \quad d_3=gcd(c,d).\]
Then the number $v$ of vertices and $e$ of egdes in a smooth decomposition is 
 \[ v = {d+2 + ( d_1+d_2+d_3) \over 2},\quad e= {3d +(d_1+d_2+d_3) \over
   2}.\]
Here, $d$ is the order of the group acting, and is arbitrary (not
necessarily odd as in \cite{RY}).
\end{lemma}
{\bf Remark:} This formula is different than that given in \cite{RY}. In
that paper, the authors only consider cases in which $d_i=1$ all $i$, which
is the same thing as only having {\it isolated} singular points, and in
these cases, our formula does agree with theirs.

\noindent{\bf Proof:} Note that if the singular point is not isolated, then
there are singular curves meeting at the point; in such a case, if the
singularity along the curve is $\fZ/e\fZ$, then $e-1$ divisors are
introduced to resolve the curve. This corresponds to $e-1$ vertices of the
decomposition, which lie on one of the edges. Let $y$ denote the number of
vertices which lie on the boundary of $\cF$, i.e., on one of the
edges. This number is determined exactly as in \cite{RY}, and is $y+3 =
d_1+d_2+d_3$, where the 3 are the original vertices of the triangle. Next,
it is easy to see that we may assume that all $y$ of these lie on one edge,
as in the following picture; the number of vertices, edges and simplices
remains contant:
\[ \epsfbox{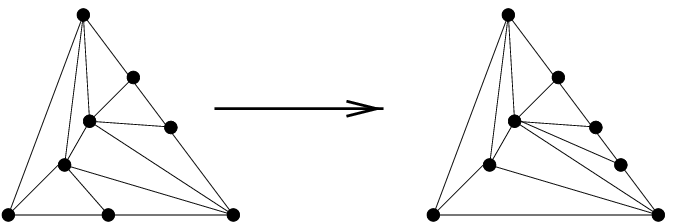}\]
Now it is easy to count the number of simplices (which is $d$) in terms of
$x=$ the number of inner vertices, and $y$. As in the following picture,
one gets the equation 
\[ x+1 + x+y = d.\]
\[ \epsfbox{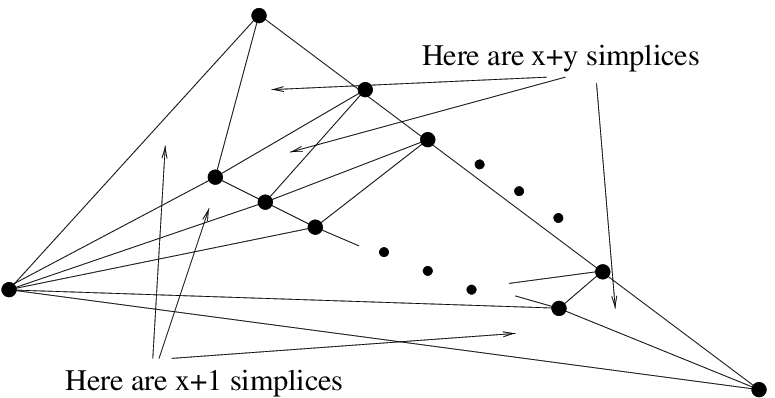} \]
Furthermore, the total number of vertices is $v=3+x+y$, and from these two
equations we get the formula for $v$. Since the triangle has Euler number =
$1=v-e+d$, the number of edges $e$ follows from this. \hfill $\Box$

\section{The twist map}
Let $V_1,\ V_2$ be weighted hypersurfaces
defined as follows.
\begin{equation}\label{8} \parbox{10cm}{$\displaystyle V_1 = \{x_0^{\ell}
    +p(x_1,\ldots,x_n) =0 \} \inn \fP_{(w_0,w_1,\ldots, w_n)} \\ V_2 = \{
    y_0^{\ell} + q(y_1,\ldots, y_m) =0 \} \inn \fP_{(v_0,v_1,\ldots,v_m)}$},
\end{equation}
where we assume both $p$ and $q$ are quasi-smooth. The degrees of these
hypersurfaces are 
\[ \nu=\deg(V_1) = \ell\cdot  w_0,\quad \mu=\deg(V_2) = \ell\cdot v_0.\]
We then consider the hypersurface 
\begin{equation}\label{9} X:= \{ p(z_1,\ldots, z_n) -q(t_1,\ldots, t_m) =0
  \} \inn \fP_{(v_0 w_1, \ldots, v_0 w_n, w_0 v_1, \ldots, w_0 v_m )}.
\end{equation}
Note that the degree of $X$ is $v_0\cdot \deg(p)=w_0\cdot \deg(q) =
v_0w_0\ell$. 
The following was shown in part I of this paper.
 The rational map 
  \begin{eqnarray*} \Phi: \fP_{(w_0,w_1,\ldots, w_n)} \times
    \fP_{(v_0,v_1,\ldots,v_m)} & \lra & \fP_{(v_0 w_1, \ldots, v_0 w_n, w_0
    v_1, \ldots, w_0 v_m )} \\ ((x_0,\ldots, x_n),(y_0,\ldots, y_m)) &
    \mapsto & (y_0^{w_1/w_0}\cdot x_1, \ldots, y_0^{w_n/w_0}\cdot x_n,
    x_0^{v_1/v_0}\cdot y_1, \ldots, x_0^{v_m/v_0}\cdot y_m)
  \end{eqnarray*}
restricts to $V_1\times V_2$ to give a rational generically finite
map onto $X$. Under the assumption that $w_0, v_0$ and $\ell$ have no
non-trivial common divisor, this map is generically $\ell$ to
one and $V_1\times V_2 \lra X$ is the projection onto the quotient of
$V_1\times V_2$ by $\boldmu_{\ell}$, which acts on the product $V_1\times
V_2$. Moreover, assuming that $X$ is Calabi-Yau, $V_2$ is Calabi-Yau and
$w_0 > 1$,
there is a resolution of singularities $\tilde{X}$ of $X$ which possesses a
fibration onto a resolution $Y$ of $V_1/\boldmu_{\ell}$ (\cite{I}, Lemma 3.4).

\section{Kodaira's 
singular elliptic fibers with torsion Monodromy}
Before turning to degenerations of K3 surfaces we show how to rederive
Kodaira's classification of singular fibers of elliptic surfaces as an
application of the twist map. So pretend we had no idea about this
classification. We will characterize singular fibers in terms of the
relations which their monodromy matrices must fulfill. In aftermath, using
the fact that the monodromy matrices are elements of $SL(2,\fZ)$, we could,
using the known properties of $SL(2,\fZ)$, derive the classification given
by Kodaira. 

These singular fibers were classified by Kodaira; his method was to
construct these fibers as quotients of smooth families of elliptic curves
of the form $D\times E$, where $D$ is a disc and $E$ is an elliptic curve
with an automorphism, i.e., of modulus either $i$ or $\grr$. His
construction was hence in terms of a {\it local} group action. We will show
how this can be easily derived upon application of our twist map, which
displays things in terms of {\it global} quotients. We
consider the following K3 surfaces which are images under the twist map of
products of a curve $C$ and an elliptic curve $E$. 

\begin{table}
\[ \!\begin{array}{|c|c|c|c|c|c|c|c|} \hline 
& (w_0,w_1,w_2) & (v_0,v_1,v_2) & \ell &  (k_1,k_2,k_3,k_4) & d & 
\hbox{singular fibers} & Monodromy  \\ \hline \hline
 1& (2,1,1) & (1,1,1) & 3 & (1,1,2,2) &  6 & 6\times \IV  & A  \\ \hline
 2&         & (1,1,2) & 4 & (1,1,2,4) &  8 & 8\times \III & B  \\ \hline
 3&         & (1,2,3) & 6 & (1,1,4,6) & 12 & 12\times \II & C  , C^2=A 
          \\ \hline
 4& (3,1,2) & (1,1,2) & 4 & (1,2,3,6) & 12 & 6\times \III, 1\times \I_0^* 
  & B,\  B^{-2}=-1  \\ \hline
 5&         & (1,2,3) & 6 & (1,2,6,9) & 18 & 9\times \II, 1\times \I_0^* & C,\ 
          C^{-3}=-1  \\ \hline
 6& (4,1,3) & (1,1,1) & 3 & (1,3,4,4) & 12 & 4\times \IV, 1\times \IV^* & A,\ 
  A^{-1}   \\ \hline
 7&         & (1,2,3) & 6 & (1,3,8,12)& 24 & 8\times \II, 1\times \IV^* & C,\ 
          C^{-2}=A^{-1} \\ \hline
 8& (5,1,4) & (1,1,2) & 4 & (1,4,5,10)& 20 & 5\times \III, 1\times \III^* & B,\
    B^{-1}  \\ \hline
 9& (7,1,6) & (1,2,3) & 6 & (1,6,14,21)&42 & 7\times \II, 1\times \II^* & C,\ 
  C^{-1}  \\ \hline
 10& (5,2,3) & (1,2,3) & 6 & (2,3,10,15)&30 & 5\times \II, 1\times \IV^*,
  1\times \I_0^* & C,\ C^{-2}, C^{-3}=-1   \\ \hline
 11& (11,5,6) & (1,2,3) & 6 & (5,6,22,33) & 66 & 2 \times \II,\ 2\times
 \II^* & C,\ C^{-1} \\ \hline  
\end{array}\]
\caption{K3 surfaces with constant modulus elliptic fibrations} 
\end{table}

In Table 1, the K3 surfaces of Fermat type with the named weights are
described as images under the twist map of products $C\times E$. There are
three elliptic curves which occur, namely:
\begin{eqnarray*} E_1 & = & \{ y_0^3 +y_1^3 +y_2^3 = 0 \} \inn
  \fP_{(1,1,1)}=\fP^2. \\ E_2 & = & \{y_0^4+y_1^4+y_2^2 = 0 \} \inn
  \fP_{(1,1,2)}. \\ E_3 & = & \{y_0^6+y_1^3+y_2^2=0\} \inn \fP_{(1,2,3)}.
\end{eqnarray*}
and the curves $C$ are those curves of
degree $w_0\ell$ in $\fP_{(w_0,w_1,w_2)}$, which we take to be of Fermat
type (except for case 11). 
In the last columns we list the Monodromy matrices, without using
what they look like. For example, in the fourth case, we have six singular
fibers of type III (this is seen by finding the number of zeros of the
polynomial $x_1^{12}+x_2^6=0\subset \fP_{(1,2)}$, which is six), and each
has monodromy matrix $B$. Recall that the monodromy gives a representation
of $\pi_1(B-\gD)$, where $B$ is the base of the fibration and $\gD$ is the
ramification locus. Since in our case $B=\fP^1$, it is known what this
fundamental group is: $\pi_1(\fP^1-\{n \hbox{ points}\}) = < \ga_1,\ldots,
\ga_n | \prod_i\ga_i =1 >$. In case 4, it follows that since we have six
singular fibers of type III, the remaining monodromy matrix (it is easily
verified that there is just one) is given by a matrix $M$ satisfying the
relation $B^6\cdot M =1$, and since $B^4=1$, it follows that
$M=B^{-2}$. But since $B^4=1$ it follows that $(B^2)^2=(B^{-2})^2=1$ and hence
$M^2=1$. Similar considerations apply in all other cases. 

To determine the structure of the degenerate fibers, note that for the first
three examples, we have (let us use $(z_1,z_2,z_3,z_4)$ as weighted
homogenous coordinates on the image weighted projective three-space) that
$z_1$ and $z_2$ are both non-vanishing, while for the sum of $d$th powers
we have $z_1^d
+z_2^d=0$. It is known that for weighted projective spaces for weights of
the form $(1,k_2,k_3,k_4)$, the affine open subset $z_1\neq 0$ is
really just a $\fC^3$. It follows that we may view, for each pair
$(z_1,z_2)$ such that $z_1^d
+z_2^d=0$, the corresponding fiber of the K3 as the curve given by the {\it
  affine} equation in $\fC^2$ 
which results. In the three cases of interest, these
affine equations are 

\[ \begin{array}{|c|c|c|} \hbox{affine equation} & \hbox{picture} &
  \hbox{Monodromy matrix}  \\   \hline  
x^3+y^2=0 & \epsfysize=1cm \epsfbox{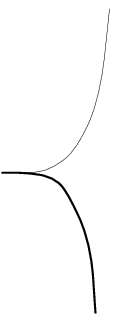} & C,\ C^6=1 \\ \hline  
x^4+y^2=0 & \epsfysize=1cm \epsfbox{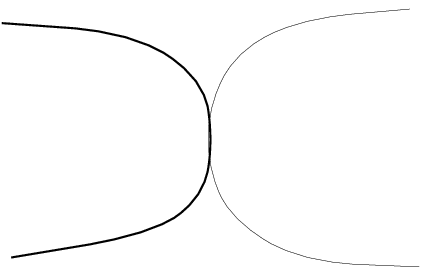} & B,\ B^4=1 \\ \hline
x^3+y^3=0 & \epsfysize=1cm \epsfbox{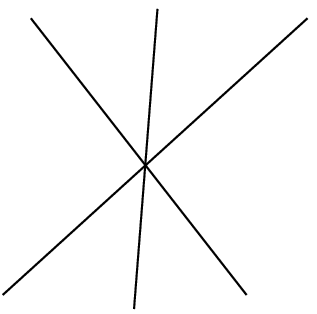} & A,\ A^3=1 \\ \hline
\end{array}
\]

To describe the other singular fibers which occur, one must consider now
the resolution of the singular weighted projective space. Then additional
singular fibers occur if: one of the ``special'' sections
$C_0:=\{z_2=0\}\cap X$ or $C_{\infty}:=\{z_1=0\}\cap X$ is {\it not} a
smooth plane cubic, an elliptic curve. 
We will describe this for cases giving rise to the
singular fibers of types I$_0^*$, II$^*$, III$^*$ and IV$^*$,
respectively. These are the cases 4 (or 5), 6, 8 and 9 above. The monodromy
matrices of these fibers are $-1,\ C^{-1},\ B^{-1}$ and $A^{-1}$,
respectively, as we have explained above. 

 {\it Case 4: $\fP_{(1,2,3,6)}$}  To see whether the curves $C_0$ and
 $C_{\infty}$ are indeed rational, we use the following compact method
 explained above: 
\[ C_0 = \fP_{(1,3,6)}[12] \cong  \fP_{(1,1,2)}[4],\] 
which is the original
elliptic curve we started with. 
\[ C_{\infty} = \fP_{(2,3,6)}[12]\cong \fP_{(2,1,2)}[4]\cong \fP_{(1,1,1)}[2],\]
which means this curve is isomorphic to a quadric in the usual projective
plane, hence {\it rational}. Next we need the singular locus of the
ambient space. This is\footnote{the notation $\{\ldots\}_{\fZ_k}$ indicates
  that $\{\ldots\}$ is fixed under a $\fZ/k\fZ$ stabilizer}
$\gS=\{z_1=z_2=0\}_{\fZ_3} \cup \{z_1=z_3=0\}_{\fZ_2} =
\gS_2\cup \gS_1$, and the equation of the hypersurface is 
$X=\{z_1^{12}+z_2^6+z_3^4+z_4^2=0\}$, 
so for the intersections we have  $\gS_1\cap X
=\fP_{(2,6)}[12]\cong \fP_{(1,3)}[6]\cong \fP_{(1,1)}[2] = 2$ points, while for
$\gS_2\cap X=\fP_{(3,6)}[12]\cong \fP_{(1,2)}[4]\cong \fP_{(1,1)}[2]=2$ points. Note
that all four points lie on the curve $C_{\infty}$, while {\it two} of them
lie on the curve $C_0$. In particular, the two curves meet in two
points. Now resolve the singularities of the ambient space; since we have
two $\fZ_2$ points and two $\fZ_3$ points, we get a total of 6 exceptional
curves, and as already mentioned, there are two ``fibers'', one rational
curve and one elliptic curve. The elliptic curve $C_0$ is clearly a smooth
fiber; it intersects two of the singular points, which has only one
explanation: there are two sections meeting it, which are components of the
exceptional locus. A picture will make this clearer:
\[ \unitlength1cm
\begin{picture}(10,7.5)
\put(-2,0){ \epsfbox{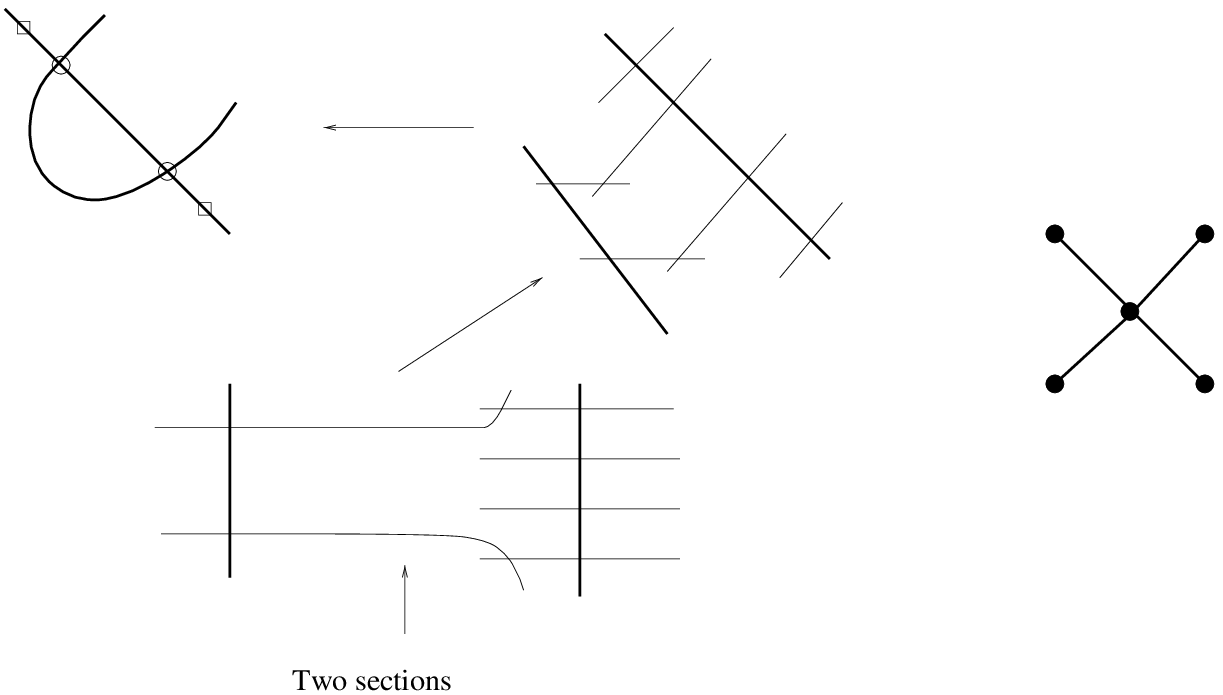} }
\put(-2.1,7.2){$C_{\infty}$}
\put(.5,6.1){$C_0$}
\put(-1,6.3){$\fZ_3$}
\put(-.2,5.7){$\fZ_3$}
\put(-2,6.4){$\fZ_2$}
\put(0,4.5){$\fZ_2$}

\put(4.1,6.8){$C_{\infty}$}
\put(4.8,3.5){$C_0$}

\put(3.5,.65){$F_{\infty} = {\bf I}_0^*$}
\put(0.3,.75){$C_0$}

\put(7,2){The extended Dynkin diagram of $D_4$}
\put(9.5,1){\large{\bf I}$_0^*$}
\end{picture}
\]

We see easily the smooth fiber $C_0$ and the fiber of type $\I_0^*$,
$F_{\infty}$, consisting of the proper transform of $C_{\infty}$ and four
exceptional $\fP^1$'s introduced in the resolution. 
Since we already deduced above that the monodromy matrix fulfills $M^2 =1$,
it follows that we have derived the structure of the singular
fiber of type I$_0^*$. 

{\it Case 6: $\fP_{(1,3,4,4)}$}
The Fermat hypersurface is $X=\{z_1^{12}+z_2^4+z_3^3+z_4^3=0\}$ and the
singular locus consists of a single component
$\gS=\{z_1=z_2=0\}_{\fZ_4}$. The intersection with $X$ is
$\fP_{(4,4)}[12]\cong \fP_{(1,1)}[3]=3$ points. Note that at these three points
the two curves $C_{\infty}$ and $C_0$ intersect. They are
$C_{\infty}=\fP_{(3,4,4)}[12]\cong \fP_{(3,1,1)}[3]$ which is rational and
$C_0=\fP_{(1,4,4)}[12]\cong \fP_{(1,1,1)}[3]$, which is the smooth elliptic
curve. There are on the intersection three $\fZ_4$ points, resolving them
gives, in addition to the three sections, one smooth fiber and one fiber of
type IV$^*$. The picture is as follows
\[\unitlength1cm
\begin{picture}(10,6.5) \put(-2,0){\epsfbox{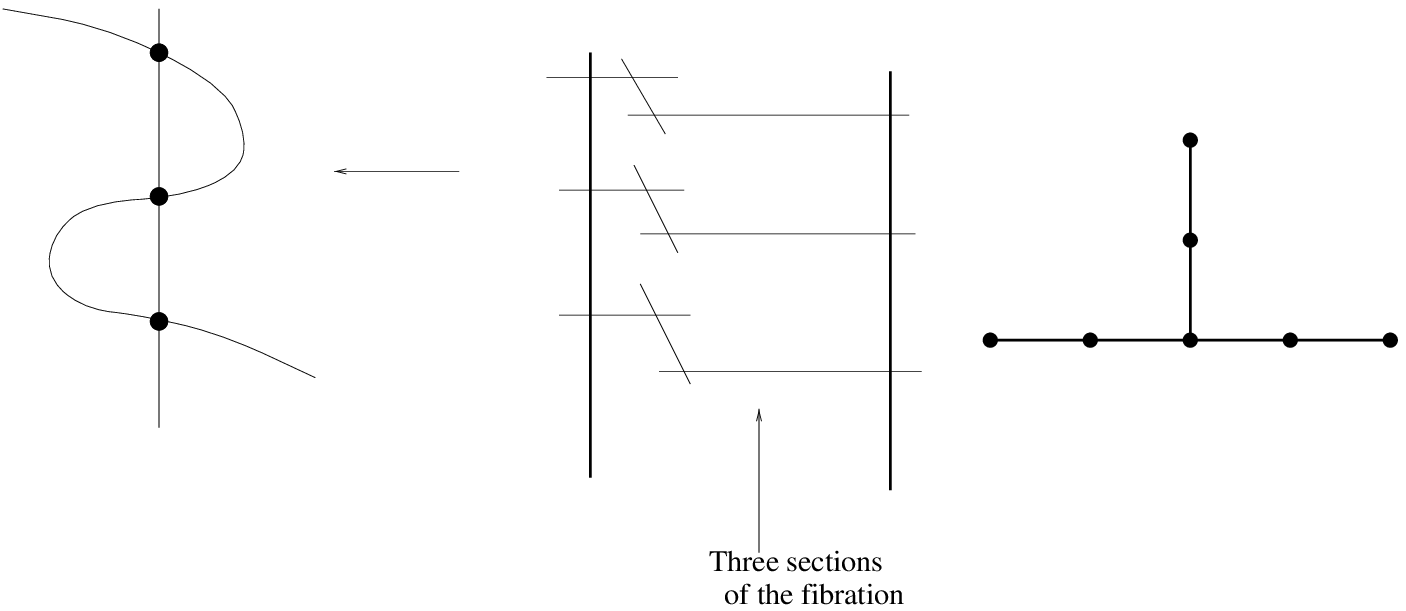}}

\put(-2.5,6){$C_{\infty}$}
\put(-.7,1.5){$C_0$}
\put(-1,5.5){$\fZ_4$}
\put(-1,4.3){$\fZ_4$}
\put(-1,3.3){$\fZ_4$}
\put(3.5,6){$C_{\infty}$}
\put(6.8,6){$C_0$}
\put(7.5,2){The extended Dynkin diagram of $E_6$}
\put(9.8,1){\large{\bf IV}$^*$}

\end{picture} 
\]

{\it Case 8: $\fP_{(1,4,5,10)}[20]$}
The singular locus here is $\gS_1=\{z_1=z_3=0\}_{\fZ_2}$ and
$\gS_2=\{z_1=z_2=0\}_{\fZ_5}$. The curve $C_0$ is $\fP_{(1,5,10)}[20]\cong 
\fP_{(1,1,2)}[4]$, which is elliptic, and $C_{\infty}=\fP_{(4,5,10)}[20]\cong 
\fP_{(4,1,2)}[4]\cong \fP_{(2,1,1)}[2]$, which is rational. The two curves $C_0$
and $C_{\infty}$ meet at $X\cap \gS_2 =\fP_{(5,10)}[20] \cong \fP_{(1,2)}[4]
\cong \fP_{(1,1)}[2]=2$ points. Hence there are two sections in the exceptional
locus. There are $X\cap \gS_1 = \fP_{(4,10)}[20] \cong \fP_{(2,5)}[10]
\cong \fP_{(1,1)}[1] =1$ point more singularities. Resolving singularities we
easily find a smooth fiber, a fiber of type III$^*$ and two sections of the
fibration. 
\[\unitlength1cm
\begin{picture}(10,3.6) \put(-3,0){\epsfbox{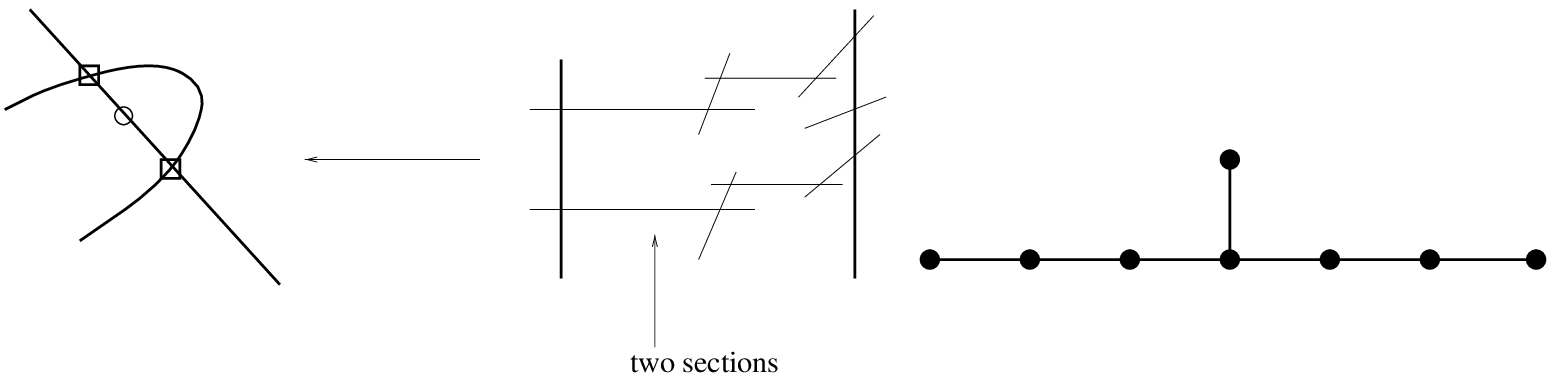}}
\put(-2.6,1.3){$C_0$}
\put(-.4,.4){$C_{\infty}$}
\put(2.5,.5){$C_0$}
\put(5.5,.5){$C_{\infty}$}
\put(-2,3.3){$\fZ_5$}
\put(-1.4,1.6){$\fZ_5$}
\put(-2.1,2.3){$\fZ_2$}

\put(6.8,.5){The extended Dynkin diagram of $E_7$}
\put(9.2,0){\large{\bf III}$^*$}

\end{picture}
\]

Finally, the most interesting case is the one giving rise to the II$^*$
type fiber.

\noindent
{\it Case 9: $\fP_{(1,6,14,21)}$} 
The Fermat hypersurface is given by 
$$X=\{ z_1^{42}+z_2^7+t_1^3+t_2^2 =0 \} \subset \fP_{(1,6,14,21)}.$$
The two special curves are given by setting $z_1$ and $z_2$ to be zero:
$$C_{\infty}=\{z_1=0\} \cap X =\{ z_2^7+t_1^3+t_2^2=0\} \subset
\fP_{(6,14,21)} = \{ (z_2')^1+(t_1')^1+(t_2')^1 =0 \} \subset
\fP_{(1,1,1)}$$
which is clearly just a linear $\fP^1$, and 
$$C_0=\{ z_1^{42} +t_1^3+t_2^2 =0\} \subset \fP_{(1,14,21)} \cong  \{ (z_1')^6
+t_1^3+t_2^2 =0\} \subset \fP_{(1,2,3)},$$
which is clearly just our elliptic curve. The singular locus of the ambient
space is $\gS_1=\{z_1=z_2=0\}_{\fZ_7}$ and $\gS_2=\{z_1=z_3=0\}_{\fZ_3}$
and $\gS_3=\{z_1=z_4=0\}_{\fZ_2}$. Furthermore, $X\cap \gS_1
=\fP_{(14,21)}[42] \cong \fP_{(2,3)}[6]\cong \fP_{(1,1)}[1]=1 $ point, $X\cap \gS_2=
\fP_{(6,21)}[42]\cong \fP_{(2,7)}[14]\cong \fP_{(1,1)}[1]=1$ point, $X\cap
\gS_3=\fP_{(6,14)}[42] \cong \fP_{(3,7)}[21]\cong \fP_{(1,1)}[1]= 1$ point. All these
points are on the curve $C_{\infty}$, one of them is the intersection with
$C_0$. This is described in the following picture

\[\unitlength1cm
\begin{picture}(10,4)
\put(-3,0){\epsfbox{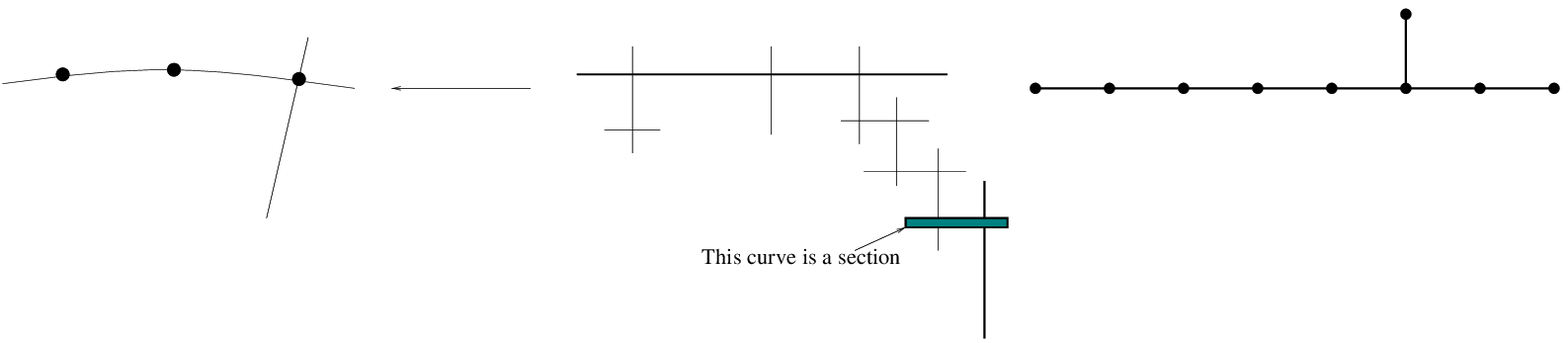} }
\put(-3,3){$C_{\infty}$}
\put(-.5,.9){$C_0$}
\put(2.2,2.85){$C_{\infty}$}
\put(6.8,-.3){$C_0$}
\put(-2.7,2.4){$\fZ_3$}
\put(-1.4,2.4){$\fZ_2$}
\put(0.2,2.3){$\fZ_7$}

\put(7.5,1.5){The extended Dynkin diagram of $E_8$}
\put(10,.5){\large{\bf II}$^*$}

\end{picture}
\]

For completeness we discuss briefly the last two cases in the table.

\smallskip
 {\it Case 10: $\fP_{(2,3,10,15)}$}. 
First we have the fixed points given by $x_0=0$. 
We are looking for solutions of 
$$ \{ x_1^{15}+x_2^{10} =0 \} \subset {\bf P}_{(2,3)},$$
which is the same as $ \{ (x_1')^5+(x_2')^5 =0 \} 
\subset {\bf P}_{(1,1)}$, of which
there are obviously only five solutions. So we have five singular fibers of
type II. The singular locus of the ambient space is
$\gS_1=\{z_1=z_2=0\}_{\fZ_5}$, $\gS_2=\{z_1=z_3=0\}_{\fZ_3}$ and
$\gS_3=\{z_2=z_4=0\}_{\fZ_2}$. We have for the intersections $X\cap \gS_1 =
\fP_{(10,15)}[30]\cong \fP_{(2,3)}[6]\cong \fP_{(1,1)}[1]=1$ point, $X\cap
\gS_2=\fP_{(3,15)}[30] \cong \fP_{(1,5)}[10]\cong \fP_{(1,1)}[2]=2$ points, and $\gS_3
= \fP_{(2,10)}[30]\cong \fP_{(1,5)}[15]\cong \fP_{(1,1)}[3]=3$ points. The two curves
$C_0$ and $C_{\infty}$ meet in a single point (the $\fZ_5$ point), hence
there is a single section in the exceptional locus. There are three $\fZ_2$
points on $C_0$, while there are two $\fZ_3$ points on the $C_{\infty}$, in
addition to the common $\fZ_5$ point. For the curves $C_0$ and $C_{\infty}$
we have $C_0=\fP_{(2,10,15)}[30]\cong \fP_{(1,5,15)}[15]\cong \fP_{(1,1,3)}[3]$,
which is a rational curve, and
$C_{\infty}=\fP_{(3,10,15)}[30]\cong \fP_{(1,10,5)}[10] \cong \fP_{(1,2,1)}[2]$,
which is also a rational curve. We have the picture:
\[\unitlength1cm
\begin{picture}(10,8)
\put(0,0){\epsfbox{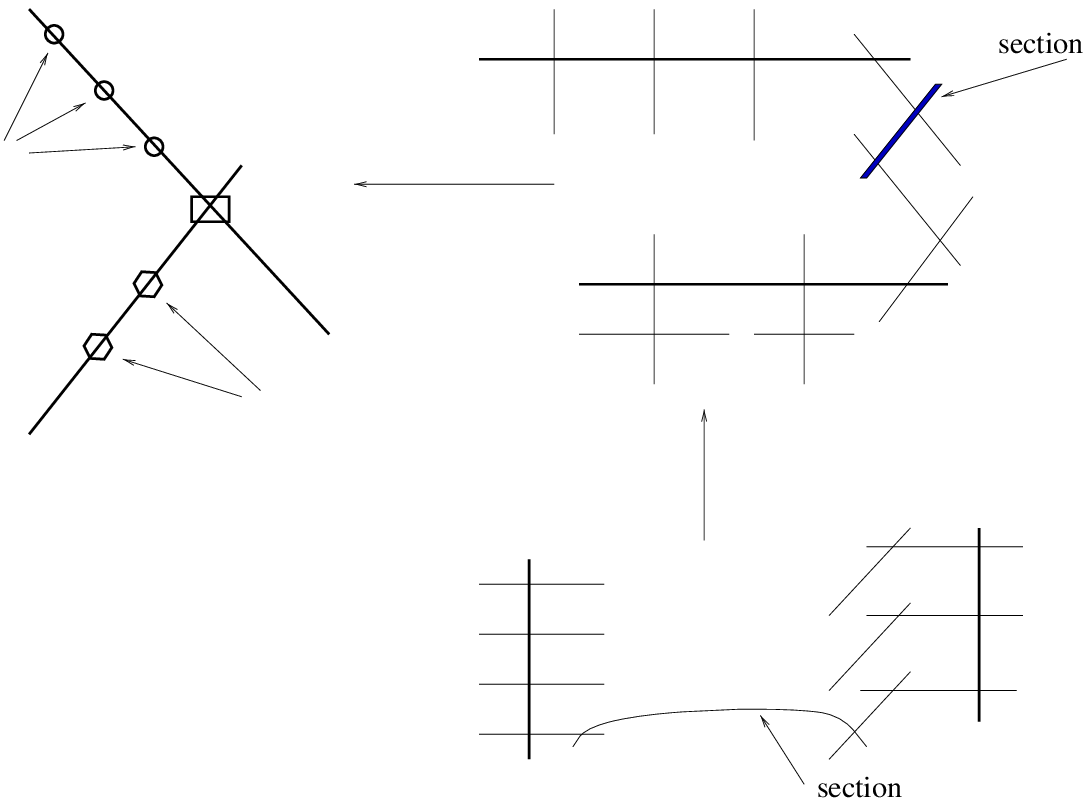}}
\put(0,3.3){$C_{\infty}$}
\put(2.9,4.5){$C_0$}
\put(4.3,7.5){$C_0$}
\put(5,5.2){$C_{\infty}$}
\put(-.8,6.6){$\fZ/2\fZ$}
\put(2.4,3.7){$\fZ/3\fZ$}
\put(.8,5.9){$\fZ/5\fZ$}
\put(5.1,2.8){$C_0$}
\put(9.8,3){$C_{\infty}$}
\put(5.1,0){I$_0^*$}
\put(9.8,.3){IV$^*$}

\end{picture}\]

\smallskip
{\it Case 11:} $\fP_{(5,6,22,33)}$
For the weights in this case, a Fermat hypersurface is
not possible. We consider instead the following polynomial:
\[  \{z_0^{12}z_1+z_1^{11}+z_2^3+z_3^2=0\} \inn \fP_{(5,6,22,33)}.\]
We see without difficulty that this is the image under the twist map
\begin{eqnarray*} \fP_{(11,5,6)}\times \fP_{(1,2,3)}  & \lra & 
 \fP_{(5,6,22,33)}\\((x_0:x_1:x_2),(y_0:y_1:y_2)) & \mapsto & 
 (y_0^{5/11}x_1: y_0^{6/11}x_2:x_0^2y_1:x_0^3y_2 ) 
\end{eqnarray*}
of the product 
$ \{ x_0^6+x_1^{12}x_2 +x_2^{11}=0\} \times \{y_0^6+y_1^3+y_2^2=0\}$.
As was explained in part I, this is a K3 surface, and we now describe the
singular fibers. The singular locus has the following components:
$\gS_1=\{z_1=z_2=0\}_{\fZ_{11}}$, $\gS_2=\{z_1=z_3=0\}_{\fZ_3}$,
$\gS_3=\{z_1=z_4=0\}_{\fZ_2}$, $\gS_4=(\gS_2\cap \gS_3)_{\fZ_6}$ and
$\gS_5=\{(1,0,0,0)\}_{\fZ_5}$. The intersections with $X$ are as
follows. $\gS_1\cap X=\fP_{(22,33)}[66]$= 1 point, $\gS_2\cap X =
\fP_{(6,33)}[66]\cong  \fP_{(6,3)}[6]\cong \fP_{(2,1)}[2] =$ 1 point, $\gS_3 \cap X =
\fP_{(6,22)}[66] \cong  \fP_{(6,2)}[6] \cong \fP_{(3,1)}[3]=$ 1 point. At the same
time the fibers $C_0$ and $C_{\infty}$ are as follows: $C_0 = \{ z_2=0\} =
\{ z_2^3+z_4^2=0\}$, which is a cusp, and $C_{\infty} =\{z_1=0\} =
\fP_{(6,22,33)}[66] \cong  \fP_{(6,2,3)}[6] \cong  \fP_{(1,1,1)}[1]$, which is a
rational curve. Note that the $\fZ_5$ point is at the cusp of $C_0$. We
have the following picture:
\[\unitlength1cm
\begin{picture}(10,8)
\put(-4,0){\epsfbox{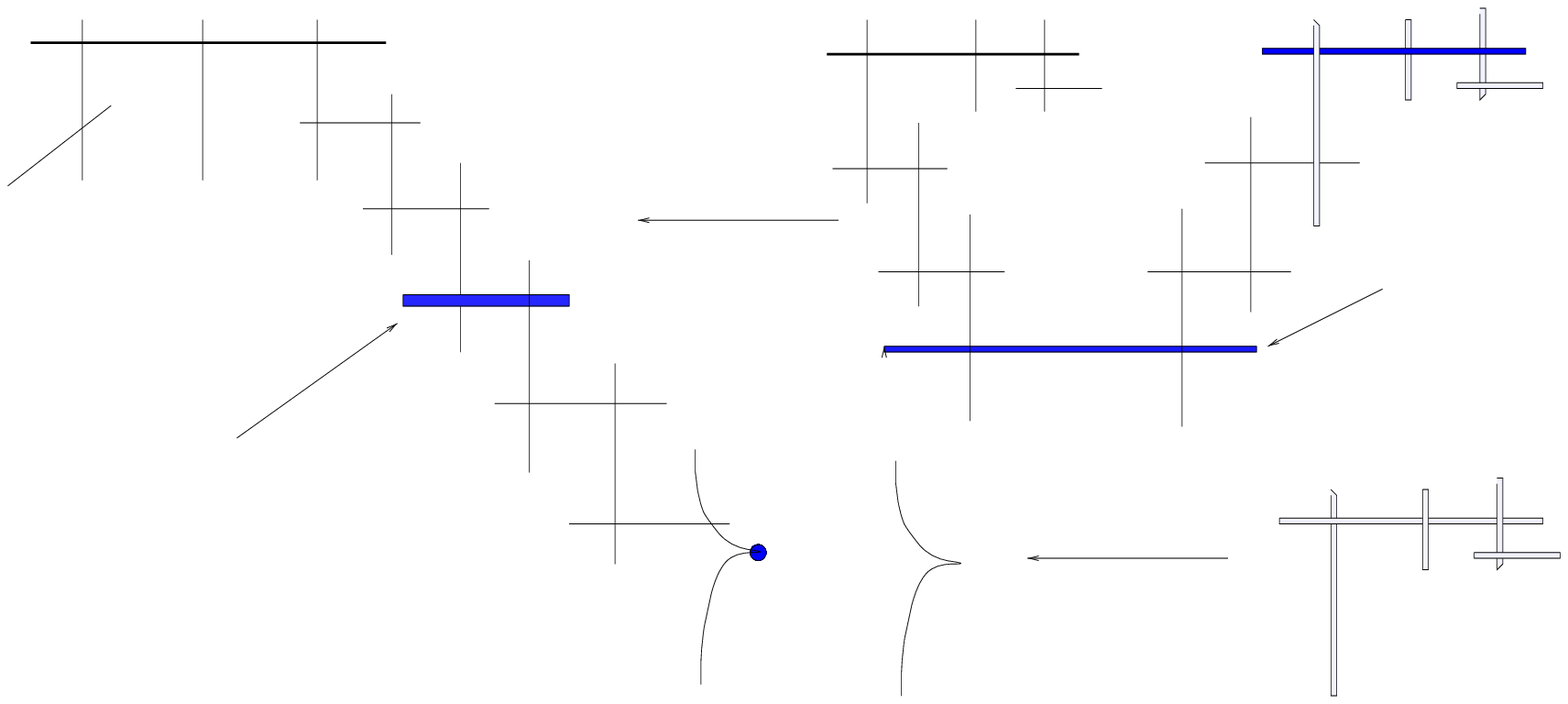}}
\put(0.6,7.2){$C_{\infty}$}
\put(3,0){$C_0$}
\put(-1.8,2.5){section}
\put(11.5,4.5){section of }
\put(11.5,4){the fibration}

\put(6.5,2.5){II$^*$}
\put(8.9,2.5){II$^*$}

\end{picture}
\]
The $\fZ/5\fZ$ point at the cusp of $C_0$ is resolved in the usual manner,
being replaced by a chain of length four. The proper transform of $C_0$ is
the ``center'' curve of the resulting configuration, which we have drawn in
the picture seperately for clarity. Thus we get the two fibers of type
$\II^*$.

\section{Singular K3 fibers with torsion monodromy}
In this section we wish to do the same as above, but now for a set of
K3-fibrations. 
\subsection{Fibers analogous to Kodaira's type II, III and IV}
First we have the analogy to the simple Kodaira fibers. Once again, we have
some unknown monodromy matrix, of which we know only the order. In these
cases, just as above, we have an affine surface as the singular
fiber. These are listed in Table 2.

\subsection{Fibers analogous to Kodaira's type I$_0^*$, II$^*$, III$^*$ and IV$^*$}

\begin{table}
\begin{tabular}{|c|c|c|c|c|c|c|} \hline 
fiber & affine equation & $\mu$ & Euler \# & fiber & Monodromy & relation \\ \hline\hline
IV$_1$ & $z^6+x^3 + y^3 =0$ & 20 & 4 & \epsfysize=2cm \epsfbox{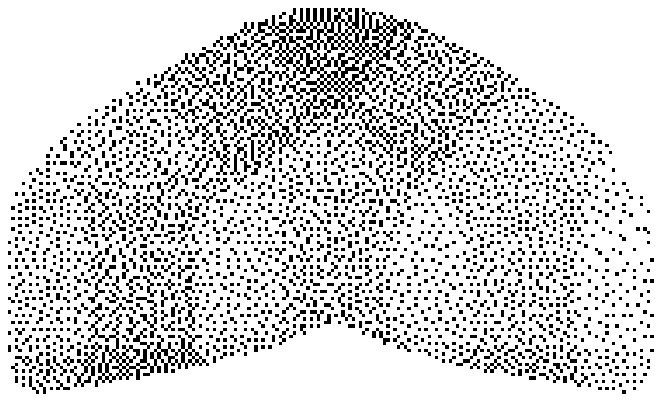} 
& M(IV$_1$) & M(IV$_1$)$^6$=1 \\ \hline
III$_1$ & $z^8+x^4+y^2=0$ & 21 & 3 & \epsfysize=2cm \epsfbox{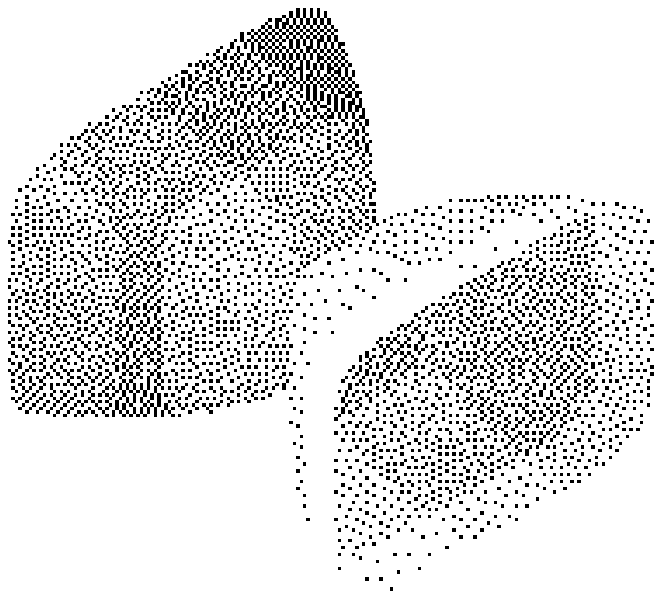} &
M(III$_1$) & M(III$_1$)$^8$=1 \\ \hline
II$_1$ & $z^{12}+x^3+y^2=0$ & 22 & 2 & \epsfysize=2cm \epsfbox{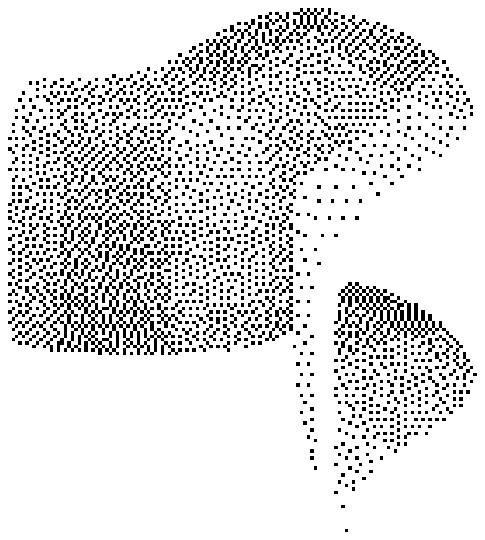} &
M(II$_1$) & M(II$_1$)$^{12}$ =1 \\ \hline
IX$_1$ & $z^6+x^4+y^2=0$ & 15 & 9 & \epsfysize=2cm \epsfbox{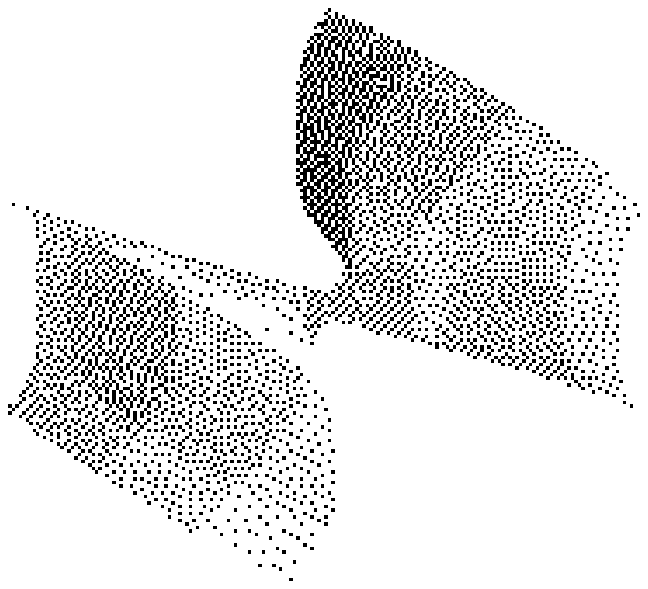} &
M(IX$_1$) & M(IX$_1$)$^{12}$=1 \\ \hline
VIII$_1$ & $z^9+x^3+y^2=0$ & 16 & 8 & \epsfysize=2cm\epsfbox{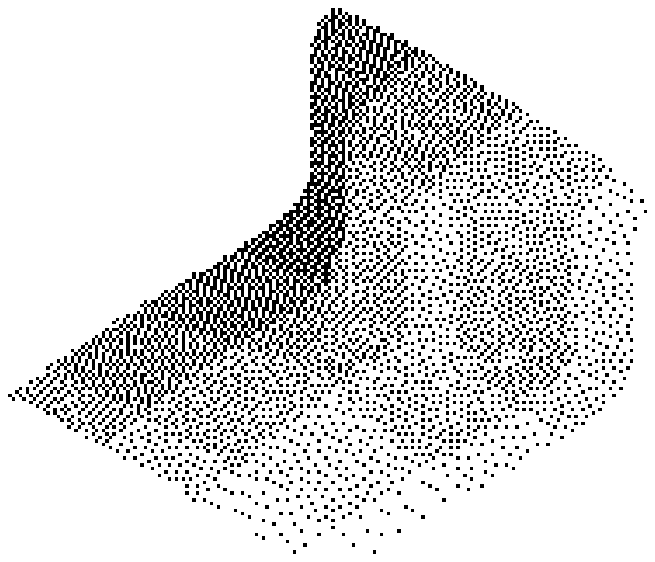} &
M(VIII$_1$) & M(VIII$_1$)$^{18}$=1  \\ \hline
XII$_1$ & $z^4+x^3+y^3=0$ & 12 & 12 & \epsfysize=2cm \epsfbox{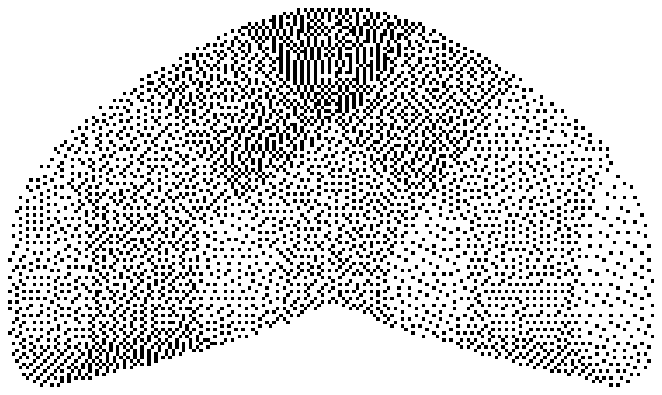} &
M(XII$_1$) & M(XII$_1$)$^{12}$=1  \\ \hline
X$_1$ & $z^8+x^3+y^2=0$ & 14 & 10 & \epsfysize=2cm \epsfbox{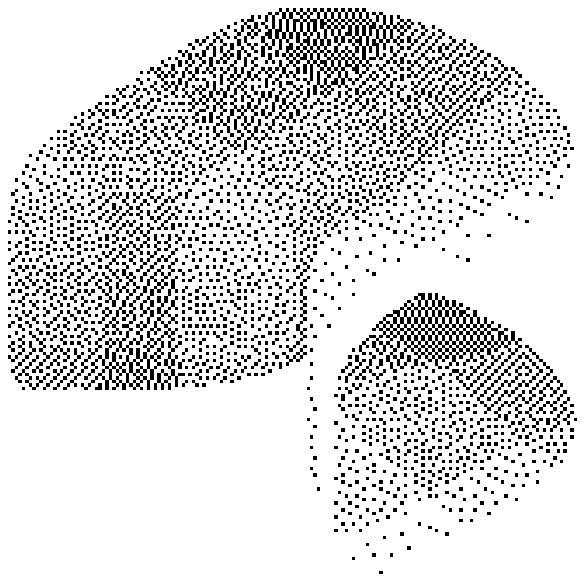} &
M(X$_1$) & M(X$_1$)$^{24}$=1 \\ \hline 
XII$_2$ & $z^5+x^4+y^2=0$ & 12 & 12 & \epsfysize=2cm \epsfbox{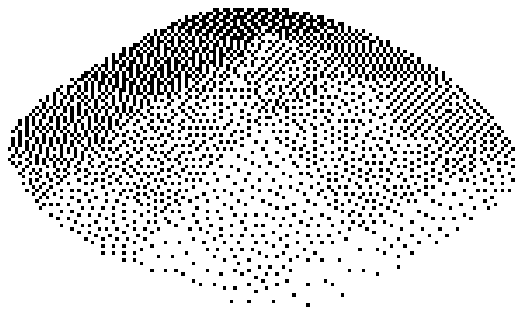}&
M(XII$_2$) & M(XII$_2$)$^{20}$=1 \\ \hline
XII$_3$ & $z^7+x^3+y^2=0$ & 12 & 12 & \epsfysize=2cm \epsfbox{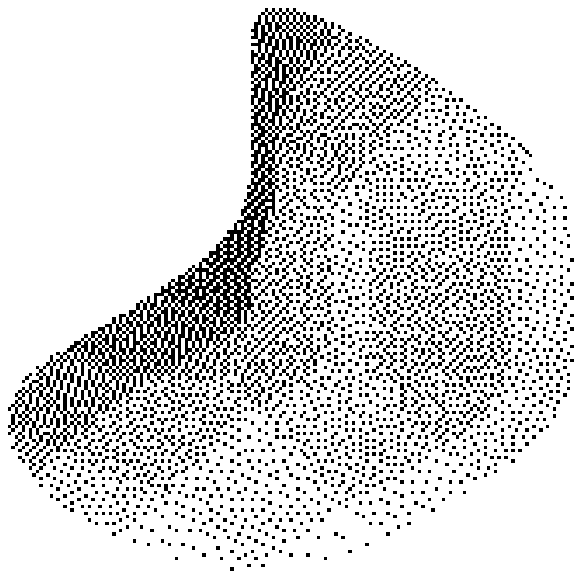} &
M(XII$_3$) & M(XII$_3$)$^{42}$=1 \\ \hline
VI$_1$ & $z^{10}+x^3+y^2=0$ & 18 & 6 & \epsfysize=2cm \epsfbox{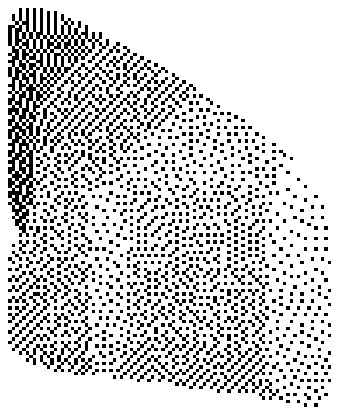}
& M(VI$_1$) & M(VI$_1$)$^{15}$=1 \\ \hline
\end{tabular}
\caption{List of singular K3-fibers with nilpotent Monodromy}
\end{table}

Here we repeat the analysis above, this time applied to weighted
hypersurfaces which are K3-fibrations. We consider the cases listed in
Table 3. 

\begin{table}
\[ \begin{array}{|l|l| c | l | c|c|c|} \hline (w_0,w_1,w_2)
    & (v_0,v_1,v_2,v_2) & \ell & (k_1,k_2,k_3,k_4,k_5) & d & Euler \# &
 singular fibers \\ \hline \hline
 (2,1,1) & (1,1,2,2) & 6 & (1,1,2,4,4) & 12 & -192 & 12\times \hbox{IV$_1$} \\
         & (1,2,3,6) & 12 & (1,1,4,6,12) & 24 & -312 & 24\times \hbox{IX}_1
         \\
         & (1,6,14,21) & 42 & (1,1,12,28,42) & 84 &  -960 & 84 
         \times \hbox{XII}_3 \\ \hline
 (3,1,2) & (1,1,2,2) & 6 & (1,2,3,6,6) & 18 & -144 & 9\times \hbox{IV}_1,
 1\times \hbox{IV}_1^*  \\
         & (1,2,3,6) & 12 & (1,2,6,9,18) & 36 & -228 & 18\times
 \hbox{IX}_1, 1\times \hbox{IX}_1^* \\
         & (1,6,14,21) & 42 & (1,2,18,42,63) & 126 & -720 & 63\times
 \hbox{XII}_3, 1\times \hbox{XII}_3^*  \\ \hline
 (4,1,3) & (1,1,2,2) & 6 & (1,3,4,8,8) & 24 & -120 & 8\times \hbox{IV}_1,\
 1\times \hbox{IV}_1^{**} \\ 
         & (1,2,3,6) & 12 & (1,3,8,12,24) & 48 & -192 & 16\times
 \hbox{IX}_1,\ 1\times \hbox{IX}_1^{**} \\
         & (1,6,14,21) & 42 & (1,3,24,56,84) & 168 & -624 & 56\times
 \hbox{XII}_3, 1\times \hbox{XII}_3^{**} \\ \hline
 (5,1,4) & (1,2,3,6) & 12 & (1,4,10,15,30) & 60 & -168 & 15\times
 \hbox{IX}_1, 1\times \hbox{IX}_1^{3*}  \\ \hline
 (7,1,6) & (1,2,3,6) & 12 & (1,6,14,21,42) & 84 & -132 & 14\times
 \hbox{IX}_1,\ 1\times \hbox{IX}_1^{4*} \\
         & (1,6,14,21) & 42 & (1,6,42,98,147) & 294 & -480 & 49\times
 \hbox{XII}_3,\ 1\times \hbox{XII}_3^{3*} \\ \hline
 (5,2,3) & (1,1,2,2) & 6 & (2,3,5,10,10) & 30 & -72 & 5\times \hbox{IV}_1,\
 1\times {\rm IV}_1^{*}, 1\times {\rm IV}_1^{**}\\
         & (1,2,3,6) & 12 & (2,3,10,15,30) & 60 & -108 &10\times
 \hbox{IX}_1, 1\times {\rm IX}_1^*, 1\times {\rm IX}^{**}_1  \\
         & (1,6,14,21) & 42 & (2,3,30,70,105) & 210 & -384 & 35\times
 \hbox{XII}_3, 1\times {\rm XII}_3^*, 1\times {\rm XII}_3^{**}  \\ \hline
       \end{array}
        \]\caption{K3-fibered Calabi-Yau weighted hypersurfaces which are
 also elliptic fibered, have constant modulus and are of Fermat type}
      \end{table}

Just as in the case of the elliptic fibrations we have the surfaces
$C_0:=\{z_2=0\}\cap X$ and $C_{\infty}:=\{z_1=0\}\cap X$. In the first
three cases it is easy to see that both of these surfaces are just smooth
fibers, hence the only singular fibers which occur are those for which the
affine surface listed in Table 2 of types IV$_1$, IX$_1$ and XII$_3$,
respectively, describe the singular fibers. Those of interest to us here
occur in the remaining cases. We begin by discussing the singular fibers
denoted IV$_1^*$, IX$_1^*$ and XII$_3^*$ in Table 3.

\subsection{Fibers analogous to I$_0^*$}
In this section we consider the cases 4-6 in the table above, in which the
monodromy matrix of the singular fiber fulfills $M^2=1$. 
In this respect, each of
these is an analog of Kodaira's I$_0^*$ type fiber.

\smallskip\noindent{\it Case }IV$_1^*$

We describe this example in more detail as a description of the general
proceedure to be used in the sequel.
The bad fiber is $C_{\infty} \cong \fP_{(2,3,6,6)}[18] \cong  \fP_{(2,1,2,2)}[6]
\cong \fP_{(1,1,1,1)}[3]$, a cubic surface. (This is a del-Pezzo surface of
degree 3, which is $\fP^2$ blown up in six points). The singular locus of
the ambient space is $\gS_1=\{z_1=z_2=0\}_{\fZ_3},\
\gS_2=\{z_1=z_3=0\}_{\fZ_2}$ and their intersection is $(\gS_1\cap
\gS_2)_{\fZ_6}$. Note that $\gS_1\cap X =\fP_{(3,6,6)}[18]\cong \fP_{(1,2,2)}[6]
\cong \fP_{(1,1,1)}[3]$, which is a cubic curve, which is elliptic. The other
intersection is $\gS_2\cap X
=\fP_{(2,6,6)}[18]\cong \fP_{(1,3,3)}[9]\cong \fP_{(1,1,1)}[3]$, which is again an
elliptic curve. They meet in the three $\fZ_6$-points on $X$
($\fP_{(6,6)}[18]\cong \fP_{(1,1)}[3] $= three points). We have the following
picture:

\[\unitlength1cm
\begin{picture}(10,3)(0,.5)\put(0,0){\epsfbox{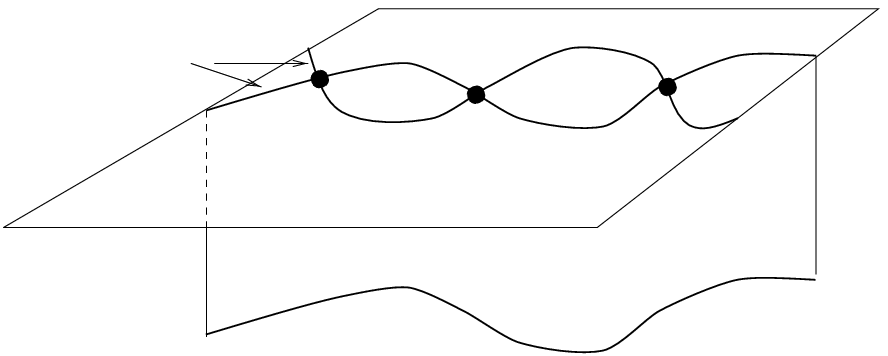}} 
\put(.9,1.5){$C_{\infty}$}\put(1.5,.5){$C_0$}\put(-.3,2.9){The curves
  $\gS_i$} 
\end{picture}
\]
Each of the dots represents a $\fZ_6$--point, $\gS_1\cap X$ is the intersection
of the two surfaces $C_0$ and $C_{\infty}$. 

The curve $\gS_1$ is the {\it base locus} of the K3-fibration, i.e., every
fiber $X_s$ passes through $\gS_1$; we think of the base $\fP^1$ of the
fibration as the exceptional $\fP^1$ of directions through $\gS_1$. Note
the general fiber is $\fP_{(1,3,6,6)}[18]\cong \fP_{(1,1,2,2)}[6]$, which
has three $\fZ_2$ points ($\{z_1=z_2=0\}_{\fZ_2} \cap X = \fP_{(2,2)}[6]
\cong \fP_{(1,1)}[3]$), whose resolutions are in fact sections of the
elliptic fibration of the fiber. This exceptional $\fP^1$ in each fiber is
the intersection of the fiber with the exceptional 
divisor ($\cong \fP_{(1,2,3)}$ described below) at each $\fZ_6$-point.

The $\fZ_6$-singular points are
of the type ${1\over 6}(1,2,3)$, which is the usual shorthand for the
quotient of $\fC_3$ by the action $(z_1,z_2,z_3)\mapsto (e^{2\pi i\over6}z_1,
e^{2\pi i \cdot 2\over 6}z_2,e^{2\pi i \cdot 3\over 6}z_3)$. 
To describe the resolution of
the singularities, we refer to the paper \cite{aspin}, in which these
$\fZ_6$-points have been resolved. The resolution is described by a cone
decomposition of the triangle with vertices $\gk_1:=(6,-2,-3),\
\gk_2:=(0,1,0),\ \gk_3:=(0,0,1)$\footnote{this is a different, but
  equivalent, description of what we described above}, 
which is the face of a three-dimensional
cone. Here this decomposition looks as follows\footnote{in this and
  following diagrams, vertices which are circled belong to divisors $E$
  which are sections of the fibration (i.e., $E\cap X_s$ is a curve for all
  $s\in \fP^1$), hence do not belong to the singular fiber.}:
\[ \epsfbox{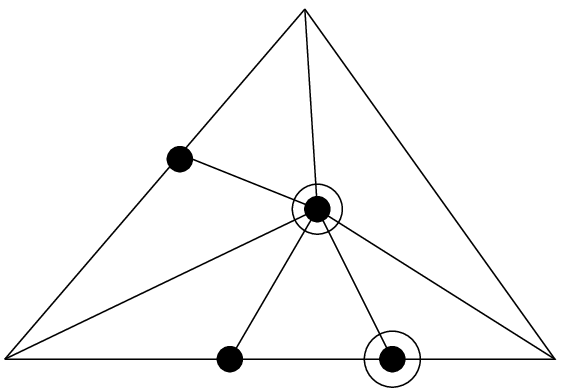} \]
in which the two (respectively one) vertices on the edge
$\overline{\gk_1,\gk_2}$ (respectively $\overline{\gk_1\gk_3}$) correspond
to the two (respectively one) exceptional divisors over $\gS_i$. One of
the two divisors over $\gS_1$ is a curve-section of the fibration: every
fiber passes through $\gS_1$, hence the intersection of one of the divisors
with each fiber is a curve, in this case isomorphic to the curve $\gS_1$,
which is elliptic. Thus, only the
other two components over the singular loci $\gS_1$ and $\gS_2$ 
belong to the singular fiber. The vertex in the middle of the
triangle corresponds to an additional exceptional divisor, which one easily
sees is just a copy of $\fP_{(1,2,3)}$, which is then resolved when the
singular curves $\gS_i$ are. There are six cones (triangles) decomposing
the big one, corresponding to the fact that the lattice is of index
six. See also \cite{RY} for details on these matters. Altogether there are
at each $\fZ_6$-point a total of four exceptional divisors; three of these
are the exceptional divisors over the $\gS_i$, the additional one at each
point is the $\fP_{(1,2,3)}$ just mentioned. Note that these latter
exceptional surfaces are also $\fP^1$-sections of the fibration, as they
lie on $\gS_1\cap X$, hence meet all fibers.

After resolution of singularities, we have the following divisors which
were introduced:
\begin{enumerate}\item Over $\gS_1$, two elliptic ruled surfaces
  $\Theta_{1,1}$ and $\Theta_{1,2}$, which intersect each other in a
  section of the ruling. \item Over $\gS_2$, an elliptic ruled surface,
  $\Theta_{2,1}$. \item Over each $\fZ_6$-point, a (resolution of a) copy of
  $\fP_{(1,2,3)}$; this intersects the fiber $F_{\infty}$ in the union of
  four rational curves, and intersects the other fibers in an exceptional
  $\fP^1$. Let $\Theta_3,\ \Theta_4,\ \Theta_5$ denote the three
  exceptional divisors introduced over the three points. \item The proper
  transform $\Theta$ of $C_{\infty}$. 
\end{enumerate}
The proper transform $[C_{\infty}]$ is the cubic surface
$\fP_{(2,3,6,6)}[18] \cong \fP_{(1,1,1,1)}[3]$ blown up at three disjoint
points. Let $\Theta$ denote this surface; it is a rational elliptic surface
with $e(\Theta)=12,\ K_{\Theta}^2=0$. The singular fiber $F_{\infty}$ is 
\[ F_{\infty} = \Theta\cup \Theta_{1,1}\cup \Theta_{2,1}.\]
The following divisors are sections of the fibration (hence do not belong
to the fiber $F_{\infty}$): $\Theta_{1,2}$, $\Theta_3$, $\Theta_4$,
and $\Theta_5$. The fiber $F_{\infty}$ is depicted in Figure 1.

\begin{figure}
\[\unitlength1cm
\begin{picture}(10,5) \put(0,0){ \epsfbox{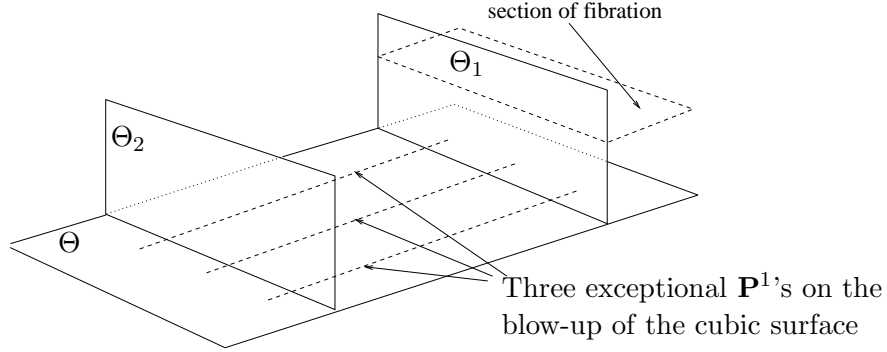}}
\put(.8,1.3){$\Theta$}\put(1.5,2.7){$\Theta_2$}\put(6,3.7){$\Theta_1$}
\put(6.7,.7){Three exceptional $\fP^1$'s on the }
\put(6.7,.2){blow-up of the cubic surface}

\end{picture}
\]
\caption{\label{figure1}This is the singular fiber of type ${\rm VI}^*_1$. 
 It consists of
  three components, the image of the K3 surface itself, here denoted
  $\Theta$, which is the proper transform of 
  a copy of $\fP_{(1,1,1,1)}[3]$, a cubic surface. The exceptional
  divisor over $\gS_i$ is denoted $\Theta_i$, and is an elliptic ruled
  surface, while $\Theta$ is a rational elliptic surface (it is $\fP^2$
  blown up at nine points), and each of the intersections
  $\Theta \cap \Theta_i$ is an elliptic curve.}
\end{figure}

We can also determine the properties of the monodromy matrix which
determines this bad fiber: since the fibration has 9 singular fibers of
type IV$_1$, and M(IV$_1$)$^6=1$, it follows that the monodromy matrix $M$
here must fulfill M(IV$_1$)$^9 \cdot M =1$, i.e., $M=$M(IV$_1$)$^{-3}$ and
hence $M^2=1$. 
So this fiber is in a sense an analoge of the Kodaira type I$_0^*$. 
Moreover, an easy calculation gives the Euler number of this fiber. Indeed,
since the fibered threefold has Euler number $-144$, and there are nine
fibers of type IV$_1$ and a single fiber of type IV$^*_1$, we get from the
formula 
\[ e(X) = 24\cdot(2-(9+1)) + 9 \cdot e({\rm IV}_1) + 1\cdot e({\rm IV}_1^*)
= -144,\]
that the Euler number of the bad fiber is $12$. We can check this, by
calculating the Euler number of our bad fiber: it is 
\[ e({\rm IV}_1^*) = e(\Theta)+e(\Theta_1)+e(\Theta_2) -e(\Theta\cap
\Theta_1) - e(\Theta\cap\Theta_2) = 12 + 0 + 0 - 0 - 0,\]
since elliptic ruled surfaces, as well as elliptic
curves, have Euler number 0. 

In the sequel we will not go into such detail.

\smallskip\noindent{\it Case} IX$_1^*$

\begin{figure}
\[\unitlength1cm
\begin{picture}(10,5)(2,0) \put(0,0){ \epsfbox{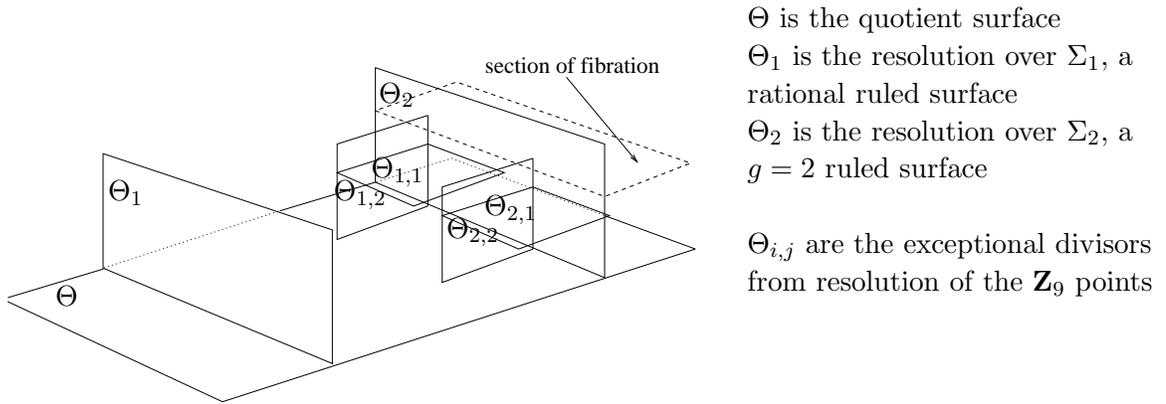}}
\put(.8,1.3){$\Theta$}\put(1.5,2.7){$\Theta_1$}\put(5.1,4){$\Theta_2$}
\put(5,3){$\Theta_{1,1}$}\put(4.5,2.7){$\Theta_{1,2}$}
\put(6.5,2.5){$\Theta_{2,1}$}\put(6,2.2){$\Theta_{2,2}$}
\put(10,5){$\Theta$ is the quotient surface}
\put(10,4.5){$\Theta_1$ is the resolution over $\gS_1$, a}
\put(10,4){rational ruled surface}
\put(10,3.5){$\Theta_2$ is the resolution over $\gS_2$, a}
\put(10,3){$g=2$ ruled surface}
\put(10,2){$\Theta_{i,j}$ are the exceptional divisors }
\put(10,1.5){from resolution of the $\fZ_9$ points}

\end{picture}\]
\caption{\label{figure2}--- The singular fiber of type IX$_1^*$ ---}
\end{figure}

The bad fiber is again $C_{\infty} \cong \fP_{(2,6,9,18)}[36] \cong 
\fP_{(2,2,3,6)}[12] \cong \fP_{(1,1,3,3)}[6]$. 
The singular locus of the ambient space is $\gS_1=\{z_1=z_2=0\}_{\fZ_3}$,
$\gS_2= \{z_1=z_4=0\}_{\fZ_2}$ and the lower-dimensional parts are given by
$\gS_3=\{z_1=z_2=z_3=0\}_{\fZ_9}$ and $\gS_4=\{z_1=z_2=x_4=0\}_{\fZ_6} =
\gS_1\cap \gS_2$. The intersections with $X$ give $\gS_1\cap X
=\fP_{(6,9,18)}[36] \cong \fP_{(2,3,6)}[12] \cong \fP_{(2,1,2)}[4]\cong \fP_{(1,1,1)}[2]$,
which is a rational curve, and $\gS_2\cap X =\fP_{(2,6,18)}[36]
\cong \fP_{(1,3,9)}[18] \cong \fP_{(1,1,3)}[6]$ which is a curve of 
genus 2 (it is a double cover of $\fP^1$ branched along a sextic, or
alternatively, a degree six curve on the Hirzebruch surface
$\fP_{(1,1,3)}$). 
These two intersect in $\gS_1\cap \gS_2\cap X =\fP_{(6,18)}[36]
\cong \fP_{(1,3)}[6] \cong \fP_{(1,1)}[2]$ which is two points. This is the locus
$\gS_4\cap X$ and consists of two points. The $\fZ_9$--locus yields
$\fP_{(9,18)}[36] \cong \fP_{(1,2)}[4] \cong \fP_{(1,1)}[2]$ which is also two
points. The configuration then looks like

\[\unitlength1cm
\begin{picture}(10,4)(0,-.5)\put(0.2,-.6){\epsfbox{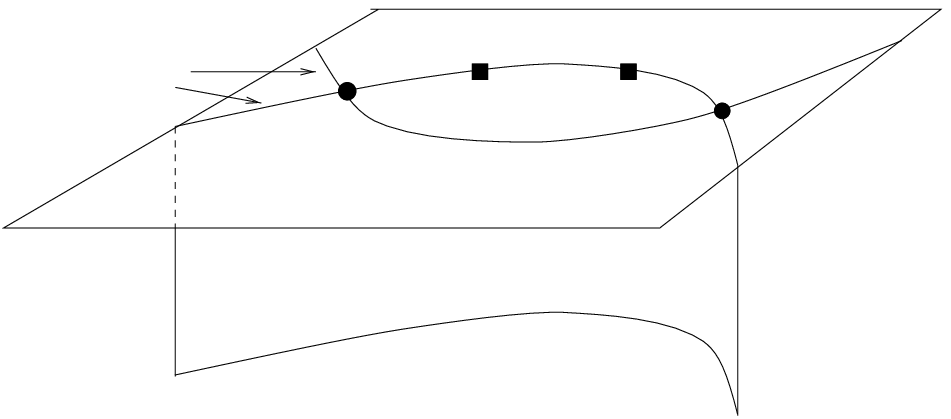}} 
\put(.9,1.5){$C_{\infty}$}\put(1.5,.5){$C_0$}\put(-.3,2.9){The curves
  $\gS_i$} 
\end{picture}
\]
in which the $\fZ_6$ points are the filled circles, the $\fZ_9$-points the
filled squares. The resolution of the curves $\gS_i,\ i=1,2$ is the same as
above, hence the resolution of the $\fZ_6$ points is precisely as above. To
describe the $\fZ_9$ points, we note that since they lie on a $\fZ_3$
curve, one of the fractions is $3/9$, hence we have ${1\over
  9}(1,3,5)$. The resolution of this follows the same pattern as above. We
now find the following cone decomposition:
\[ \epsfbox{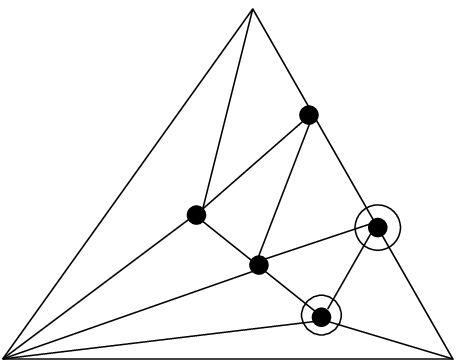} \]
Using Lemma \ref{formula}, 
we see that in our case of ${1\over 9}(1,3,5)$ we have $d_1=d_3=1,\
d_2=3,\ d_1+d_2+d_3=5$ and hence the number of vertices (including the
three corners of the original triangle) is 8, the number of edges is 16,
and the cone decomposition is as given above. The two vertices on the edge
of the triangle of course correspond to the two exceptional divisors over
the curve $\gS_1$, on which the $\fZ_9$-points lie. There are at each
$\fZ_9$-point three additional exceptional divisors. Note that again, one
of these three is a $\fP^1$-section of the fibration, while the other two
are components of the bad fiber. We have circled the vertices corresponding
to the exceptional divisors which are $\fP^1$-sections of the fibration. 
After the resolution of singularities, our fiber will look as in Figure
\ref{figure2}. 

Once again the monodromy is easily seen to satisfy $M^2=1$, 
and the Euler number of
this singular fiber can be calculated as above; it is 18.

\smallskip\noindent{\it Case} XII$_3^*$

The ambient space is $\fP_{(1,2,18,42,63)}$, and the singular locus in this
space is 
\begin{itemize} \item $\gS_1=\{z_1=z_2=0\}_{\fZ_3} \cong \fP_{(18,42,63)}
  \cong \fP_{(1,1,3)}$. \item $\gS_2 =\{ z_1=z_5=0\}_{\fZ_2} \cong
  \fP_{(2,18,42)} \cong \fP_{(1,3,7)}$. \item $\gS_3=\gS_1\cap \gS_2
  =\{z_1=z_2=z_5=0\}_{\fZ_6} \cong \fP^1$. 
  \item $\gS_4=\{z_1=z_2=z_3=0\}_{\fZ_{21}} \cong \fP^1$. 
  \end{itemize}
The bad fiber is $C_{\infty}\cong \fP_{(2,18,42,63)}[126] \cong
\fP_{(1,3,7,21)}[21]$; it contains the two curves $X\cap \gS_1$ and $X\cap 
\gS_2$, and these two curves intersect in just one point. There is a
further singular point $\gS_4\cap X$, which also lies on $\gS_1$ but not on
$\gS_2$. Hence we have the picture
\[\unitlength1cm
\begin{picture}(10,4)(0,.5)\put(0,0){\epsfbox{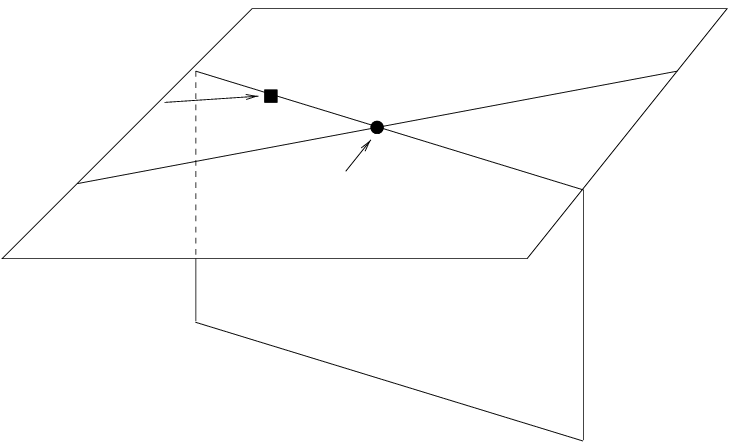}}
\put(3,4){$C_{\infty}$}\put(3,1){$C_0$}
\put(7,3.7){$\gS_2$, a $\fZ_2$ curve}
\put(6.2,2.5){$\gS_1$, a $\fZ_3$ curve}
\put(2.2,2.5){The $\fZ_6$ singular point} 
\put(-2,3.4){The $\fZ_{21}$ singular point}

\end{picture}
\]
The $\fZ_6$-point is resolved just as above, with a single exceptional
divisor which is a $\fP^1$-section of the fibration. For the $\fZ_{21}$
point, we deduce from the fact that it lies on a $\fZ_3$ curve, that it is
of type ${1\over 21}(1,2,18)$. This yields $d_1=d_2=1,\ d_3=3$ and hence,
by the formula above, $v={23+5 \over 2} = 14$, of which five are on the
boundary of the triangle. As we have already mentioned, any decomposition
of the cone with this many vertices yields a smooth resolution of the
singular point. As described above, since the singular point
is on $\gS_1$, it meets every fiber, so one of the components is a
curve-section of the fibration. We choose the following decomposition:
\[ \epsfbox{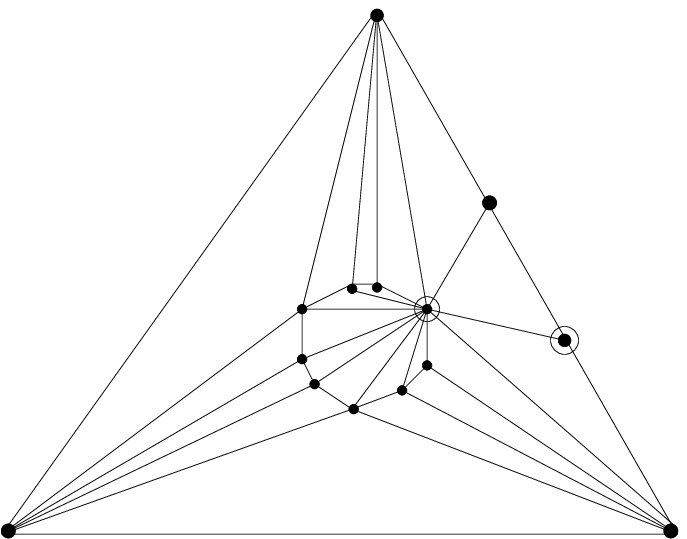} \]

\begin{figure}
\[\unitlength1cm \begin{picture}(10,4.3)(2,-.2)
\put(0,0){\epsfbox{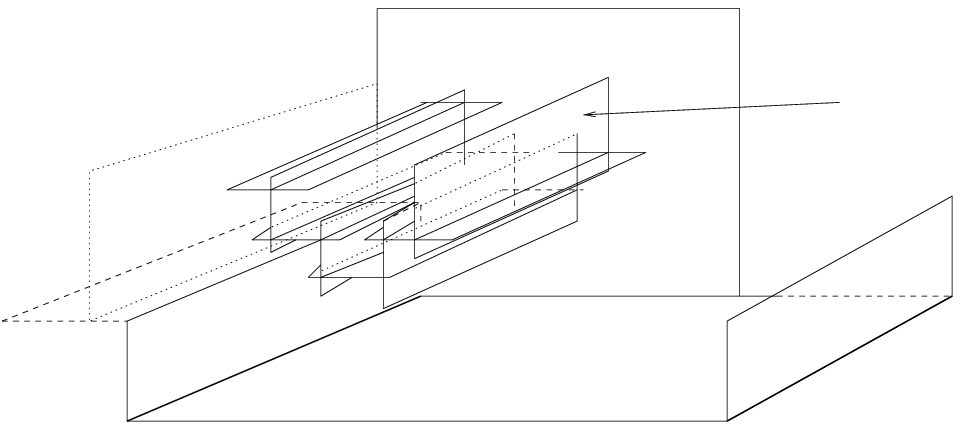}}
\put(1.4,.65){$\Theta_{1,1}$}\put(.5,.75){$\Theta_{1,2}$}
\put(2.3,.2){$\Theta$} \put(.7,3){$C_0$} \put(5,3.6){$\Theta_9$}
\put(8,3.5){$\Theta_1-\Theta_8$}
\put(7.5,.5){$\Theta_{2,1}$}
\put(10,4){The fiber is the sum }
\put(10,3.5){$\Theta + \Theta_1+\ldots,
  +\Theta_8+\Theta_{1,1}+\Theta_{2,1}$}
\put(10,3){where $\Theta$ is the quotient fiber}
\put(10,2.5){ $\cong
  \fP_{(1,3,7,21)}[21]$, $\Theta_i,\ i=1,\ldots, 8$ }
\put(10,2){are the resolution divisors of }
\put(10,1.5){the $\fZ_{21}$ point, $\Theta_9$ the one which}
\put(10,1){is a section, $\Theta_{1,1}$ is rational ruled}
\put(10,.5){$\Theta_{2,1}$ is $g=6$ ruled.}

\end{picture}
\]
\caption{\label{figure3} --- The singular fiber of type XII$^*_3$ ---}
\end{figure}

We have nine exceptional divisors, of which one is a section of the
fibration, the other eight belong to the singular fiber. After resolution
of singularities, our singular fiber thus looks as in Figure \ref{figure3}.

Note the surface $\Theta_9$, which is a section of the fibration, not a
component of the singular fiber. This component of the resolution
corresponds to the circled vertex above, and meets all the other eight
components, while each of the other eight meets only two others, as drawn.

\subsection{Fibers analogous to Kodaira's type IV$^*$}
In this section we consider the cases 7-9 in our table. In all these cases,
the monodromy matrix is an element of order 3, so in a sense an analog of
Kodaira's type IV$^*$ type fiber.

\begin{figure}
\[ \unitlength1cm \begin{picture}(10,4)(2,0)
\put(-.5,0){
\put(0,0){\epsfbox{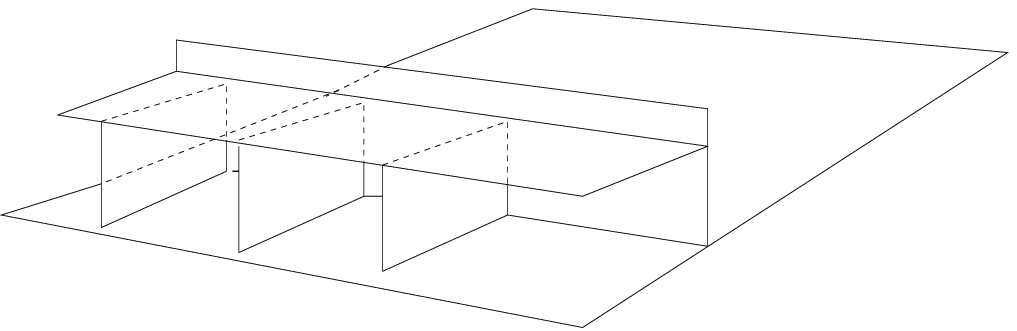}}
\put(8,2.5){$\Theta$}\put(6.5,1.2){$\Theta_{1,1}$}
\put(5.3,1.7){$\Theta_{1,2}$} \put(1.2,1.5){$\Theta_1$}
\put(2.8,1.3){$\Theta_2$} \put(4.3,1.2){$\Theta_3$}
}
\put(10,4){The fiber is the sum }
\put(10,3.5){$\Theta+\Theta_1+\Theta_2+\Theta_3+\Theta_{11}+\Theta_{1,2}$}
\put(10,3){where $\Theta$ is the quotient fiber }
\put(10,2.5){$\cong \fP_{(1,2,2,3)}[6]$, $\Theta_i,\ i=1,\ldots, 3$}
\put(10,2){are the resolution divisors of the }
\put(10,1.5){three $\fZ_8$ points, $\Theta_{1,i},\ i=1,2$}
\put(10,1){are the resolution divisors}
\put(10,.5){over the singular curve.}

\end{picture}
\]
\caption{\label{figure4} --- The singular fiber of type IV$^{**}_1$ ---}
\end{figure}

\smallskip\noindent{\it Case } IV$^{**}_1$:

The ambient space is $\fP_{(1,3,4,8,8)}$, which has the following
singularities:
\begin{itemize}\item $\gS_1 = \{z_1=z_4=0\}_{\fZ_4}$, \item
  $\gS_2=\{z_1=z_2=z_3=0\}_{\fZ_8}$. 
\end{itemize}
For the intersections we have $\gS_1\cap X = \fP_{(4,8,8)}[24] \cong
\fP_{(1,2,2)}[6]\cong \fP_{(1,1,1)}[3]$, which is an elliptic curve. $\gS_2
\cap X = \fP_{(8,8)}[24]\cong \fP_{(1,1)}[3]$, which consists of three
points. So the singular locus of $X$ is a curve and three additional points
on that curve. The bad fiber is $C_{\infty} = \fP_{(3,4,8,8)} \cong
\fP_{(1,2,3,2)}[6]$, which is a double cover of the space $\fP_{(1,2,3)}$
branched over a sextic curve. Since the $\fZ_8$ points lie on a $\fZ_4$
curve, their type is ${1\over 8}(1,3,4)$, and we have $d_1=d_2=1,\ d_3=4$
so by our formula above we get $v={10+6 \over 2} =8$, of which 3 lie on the
edge (corresponding to the exceptional divisors of the $\fZ_4$-curve), 
and we take the following cone decomposition:
\[ \epsfbox{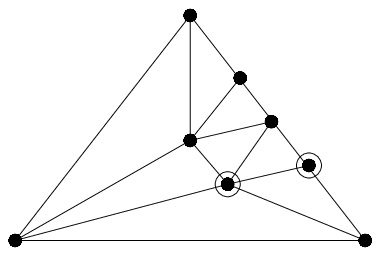} \]
This means that at each of the three $\fZ_8$ points, we have two
exceptional divisors, one of which is a section of the fibration. In the
picture of the cone decomposition we have circled the two vertices which
correspond to exceptional divisors which are sections of the fibration. The
others then belong to the singular fiber, so we have a total of five
components in addition to the image of the K3 itself. After resolution of
singularities, we get the fiber in Figure \ref{figure4}.

\begin{figure}
\[\unitlength1cm\begin{picture}(10,4)(.7,0)
 \put(-2,-.6){\epsfbox{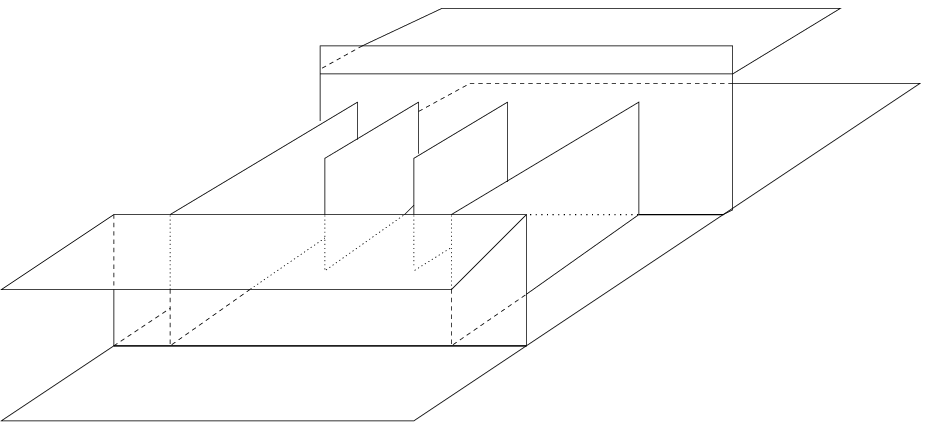}}
\put(6,2.5){$\Theta$}\put(4.7,2){$\Theta_{2,1}$}
\put(5,3.4){$\Theta_{2,2}$}\put(4,1.8){$\Theta_1$}
\put(1.5,1.8){$\Theta_3$}\put(.7,1.8){$\Theta_2$}
\put(2.5,1.9){$\Theta_4$}
\put(2.7,.6){$\Theta_{1,1}$}
\put(-1,1){$\Theta_{1,2}$}
\put(5.7,1.5){$\gS_1$}\put(3.7,.1){$\gS_2$}

\put(8,4){The singular fiber is the sum}
\put(8,3.5){$\Theta+\Theta_1+\Theta_2+\Theta_3+\Theta_4+$}
\put(8,3){$\Theta_{1,1}+\Theta_{1,2}
  +\Theta_{2,1}+\Theta_{2,2}$}
\put(8,2.5){where $\Theta$ is the proper transform of}
\put(8,2){the surface $C_{\infty}$ and the $\Theta_i$}
\put(8,1.5){$i=1,\ldots,4$ are the resolution surfaces of }
\put(8,1){the two $\fZ_8$ points and the }
\put(8,.5){the two $\fZ_{12}$ points; the $\Theta_{1,i},\ i=1,2$}
\put(8,0){are rational ruled, the $\Theta_{2,i},\ i=1,2$}
\put(8,-.5){are elliptic ruled.}

\end{picture}
\]
\caption{\label{figure5}--- The singular fiber of type IX$^{**}_1$ ---}
\end{figure}

\bigskip\noindent{\it Case} IX$^{**}_1$:

The ambient space is $\fP_{(1,3,8,12,24)}$, which has as singular locus:
\begin{itemize}\item $\gS_1=\{z_1=z_2=0\}_{\fZ_4}$, \item
  $\gS_2=\{z_1=z_3=0\}_{\fZ_3}$, \item $\gS_3=\gS_1\cap \gS_2
  =\{z_1=z_2=z_3=0\}_{\fZ_{12}}$, \item $\gS_4=\{z_1=z_2=z_4=0\}_{\fZ_8}$. 
\end{itemize}
The intersections with $X$ are $\gS_1\cap X = \fP_{(8,12,24)}[48]\cong
\fP_{(1,1,1)}[2]$, a rational curve, $\gS_2\cap X =\fP_{(3,12,24)}[48]\cong
\fP_{(1,1,2)}[4]$, an elliptic curve, $\gS_3\cap X = \fP_{(12,24)}[48]\cong
\fP_{(1,1)}[2]$, consisting of two points, and $\gS_4\cap X =
\fP_{(8,24)}[48]\cong \fP_{(1,1)}[2]$, again two points. Hence, before
resolution, our fiber $C_{\infty} = \fP_{(3,8,12,24)}[48] \cong
\fP_{(1,1,2,2)}[3]$ looks as follows:
\[ \unitlength1cm\begin{picture}(10,3.5)
\put(0,0){\epsfbox{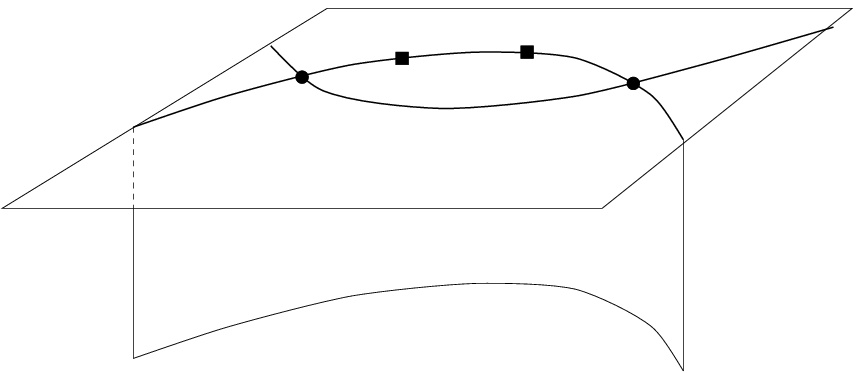}}
\put(2,1){$C_0$}\put(2.5,2){$C_{\infty}$}
\put(0,2.5){$\gS_1 : \fZ_4$}
\put(1.1,3.3){$\gS_2 : \fZ_3$}
\put(2.8,2.7){$\fZ_{12}$}\put(6.2,2.6){$\fZ_{12}$}
\put(3.9,3.4){$\fZ_8$}\put(5.3,3.4){$\fZ_8$}

\end{picture}
 \]
The situation at each of the $\fZ_8$ points is just as in the previous
example, as each lies again on a $\fZ_4$ curve. Hence the resolution of
these points introduces two exceptional divisors each, one of which is a
component of the singular fiber. Looking at the $\fZ_{12}$ points, we see
that as they are again on $\fZ_3$ curves, they are of type ${1\over
12}(1,2,9)$. Hence we have $d_1=1,\ d_2=2,\ d_3=3$ and the formula for the
number of vertices of a cone decomposition yields $v= {14+6\over 2}=10$
vertices, of which there are three on the boundary corresponding to the $\fZ_4$
curve, and two others on another boundary, corresponding to the $\fZ_3$
curve. Hence there are two vertices in the interior, and we have the cone
decomposition:
\[\epsfbox{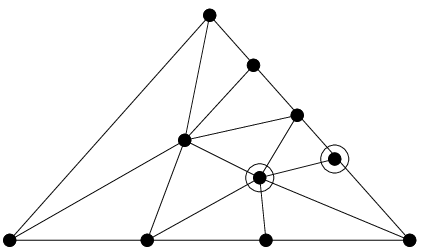}\]
In total, in addition to the proper transform $\Theta$ of the singular
fiber $C_{\infty}$, we have two divisors over $\gS_1$ (the third is a
section of the fibration and not a component of the fiber), two over
$\gS_2$, and over each of the two $\fZ_8$ and two $\fZ_{12}$ points we also
have two exceptional divisors, one of which is not a section, hence a
component of the fiber. After resolution of singularities, we have the
picture of Figure \ref{figure5}.

\bigskip\noindent{\it Case} XII$^{**}_3$:

The ambient space is $\fP_{(1,3,24,56,84)}$, which has the following
singular locus:
\begin{itemize}\item $\gS_1=\{z_1=z_2=0\}_{\fZ_4}$, \item
  $\gS_2=\{z_1=z_4=0\}_{\fZ_3}$, \item $\gS_3=\{z_1=z_2=z_5=0\}_{\fZ_8}$,
  \item $\gS_4=\{z_1=z_2=z_3=0\}_{\fZ_{28}}$, \item $\gS_5=\gS_1\cap \gS_2
      =\{z_1=z_2=z_4=0\}_{\fZ_{12}}$. 
\end{itemize}
For the intersections with $X$ we get $\gS_1\cap X\cong\fP_{(1,1,1)}[1]$ is
rational, $\gS_2\cap X\cong \fP_{(1,2,7)}[14]$, a $g=3$ curve, $\gS_3\cap
X$, $\gS_4\cap X$ and $\gS_5\cap X$ all consist of just a single point. The
singular fiber is $C_{\infty}\cong \fP_{(1,2,7,14)}[14]$, a rational
surface. This looks as follows:

\[ \unitlength1cm \begin{picture}(10,4)
  \put(0,0){\epsfbox{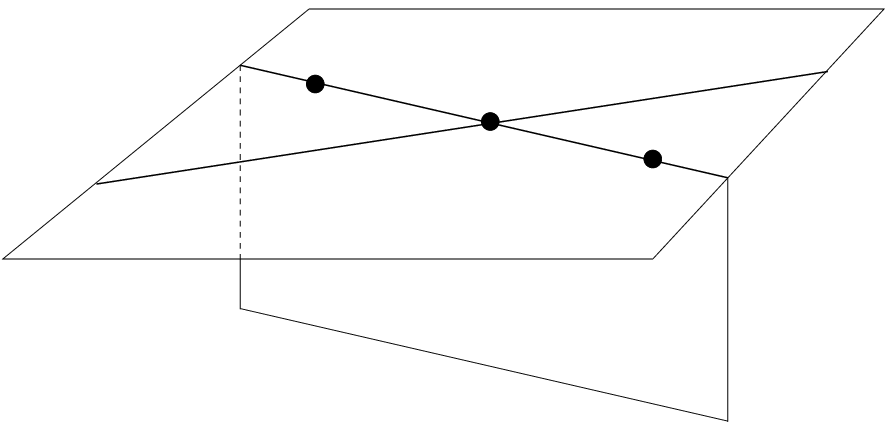}}
\put(6,3.9){$C_{\infty}$}\put(6.5,.5){$C_0$}
\put(3.7,3.6){$\fZ_8$}\put(6.5,2.3){$\fZ_{28}$}
\put(5,2.6){$\fZ_{12}$} \put(7.5,2.3){$\gS_1 : \fZ_4$}
\put(8.5,3.3){$\gS_2 : \fZ_3$}

\end{picture}
\]

\begin{figure}
\[\unitlength1cm \begin{picture}(10,5)(3.5,0) 
\put(0,0){\epsfbox{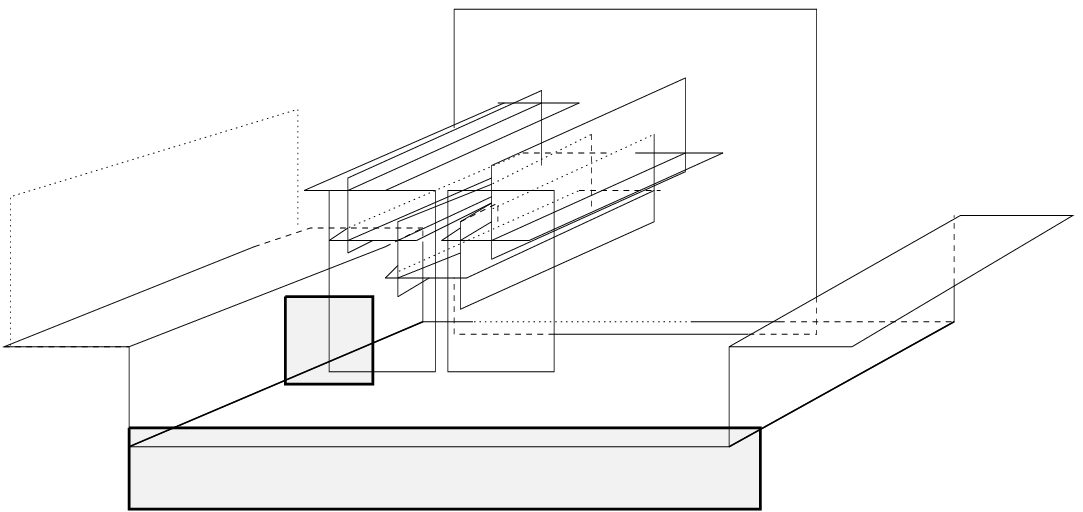}}
\put(2.7,1){$\Theta$}\put(1.5,1.2){$\Theta_{1,1}$}
\put(.9,1.8){$\Theta_{1,2}$}\put(7.5,1.3){$\Theta_{2,1}$}
\put(8.6,2.2){$\Theta_{2,2}$} \put(7.2,4.2){$\Theta_9$}
\put(3.7,2){$\Theta_{10}$}\put(5,2){$\Theta_{11}$}
\put(2.7,.25){$\Theta_{12}$}\put(3.2,1.5){$\Theta_{13}$}

\put(10,4){The singular fiber is the sum}
\put(10,3.5){$\Theta+\Theta_1+\cdots+\Theta_{13}+\Theta_{1,1}+ \Theta_{1,2}
  +\Theta_{2,1} +\Theta_{2,2}$}
\put(10,3){where $\Theta$ is the proper transform of $C_{\infty}$}
\put(10,2.5){and $\Theta_i,\ i=1,\ldots, 11$ are the divisors}
\put(10,2){of the resolution of the $\fZ_{28}$ point}
\put(10,1.5){$\Theta_{12}$ is from the $\fZ_{12}$ point}
\put(10,1){$\Theta_{13}$ is from the $\fZ_{8}$ point}
\put(10,.5){$\Theta_{1,i},\ i=1,2$ are rational ruled}
\put(10,0){$\Theta_{2,i}\ i=1,2$ are $g=3$ ruled. }

\end{picture}
\]
\caption{\label{figure6} --- The singular fiber of type XII$^{**}_3$ --- }
\end{figure}

The $\fZ_{12}$ point is just the same as above, yielding upon resolution
one additional component to the fiber. Similarly, the $\fZ_8$ point is the
same as above, yielding also one additional component to the fiber. It
remains to resolve the $\fZ_{28}$ point. Note that as it lies on a $\fZ_4$
curve, it is of type ${1\over 28}(1,3,24)$. We have $d_1=d_2=1,\ d_3=4$ and
for the number of vertices we get $v={30+6\over 2} =18$, of which only six
are on the boundary. Hence we must insert 12 vertices in the interior. We
choose the following cone decomposition:

\[ \epsfbox{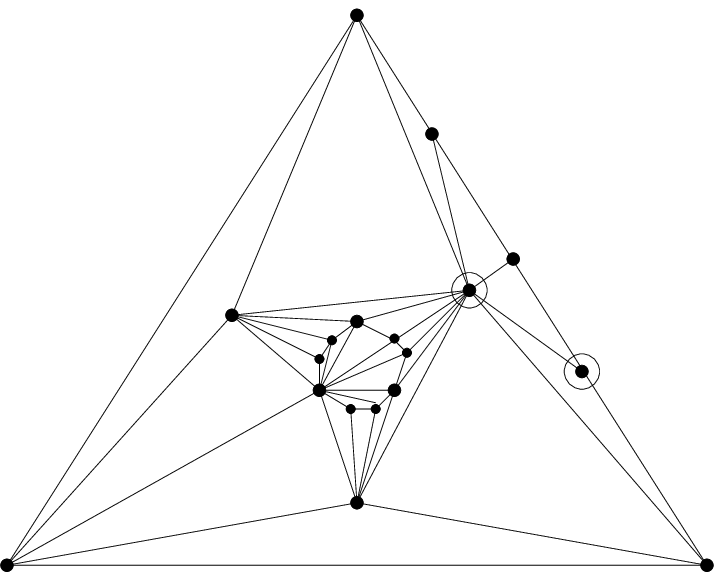} \]
Once again, the two components of the resolution which give rise to
sections of the fibration instead of components of the singular fiber are
circled. We again see in the middle the eight components giving rise to an
$A_7$ configuration, but this time the component which meets all
eight is a component of the singular fiber instead of a section. In
addition, we have two more components, each of which meets the special
component and four of the components of the $A_7$ chain. After
resolution, the fiber looks as in Figure \ref{figure6}.

\subsection{An analog of Kodaira's type III$^*$ fiber}
We consider now the tenth case in the table. The monodromy matrix has order
four, so in this sense this is an analog of Kodaira's type III$^*$ fiber. 

\smallskip\noindent{\it Case} IX$^{***}_1$:

The ambient space is $\fP_{(1,4,10,15,30)}$ with singular locus
\begin{itemize}\item $\gS_1=\{z_1=z_2=0\}_{\fZ_5}$, \item
  $\gS_2=\{z_1=z_4=0\}_{\fZ_2}$, \item $\gS_3=\gS_1\cap \gS_2 =
  \{z_1=z_2=z_4=0\}_{\fZ_{10}}$, \item $\gS_4 = \{
  z_1=z_2=z_3=0\}_{\fZ_{15}}$. 
\end{itemize}
For the intersections with $X$ we have $\gS_1\cap X \cong
\fP_{(1,1,1)}[2]$, a rational curve, $\gS_2\cap X \cong \fP_{(1,2,3)}[6]$,
an elliptic curve, and the two intersections $X\cap \gS_3$ and $X\cap
\gS_4$ both consist of two points. Hence the picture is just as in the case
IX$^{**}_1$ above, with the $\fZ_{12}$ points now replaced by $\fZ_{10}$
points, and the $\fZ_{8}$ points there replaced by $\fZ_{15}$ points
here. The resolution of the $\fZ_2$ curve is an elliptic ruled surface, a
component of the singular fiber, and the resolution of the $\fZ_5$ curve is
a union of four rational ruled surfaces, of which one is a section of the
fibration, while the other three give components of the singular fiber. So
we just have to resolve the $\fZ_{10}$ and $\fZ_{15}$ points. Note that
since both lie on a $\fZ_5$ curve, they are of types ${1\over 10}(1,4,5)$
and ${1\over 15}(1,4,10)$ (or $(2,3,x)$, it doesn't matter). Hence we have
$d_1=1,\ d_2=2,\ d_3=5$ in the first and $d_1=d_2=1,\ d_3=5$ in the second
case. Hence we have $v={12+8\over 2}=10$ vertices in the first case and
$v={17+7 \over 2}=12$ vertices in the second. We take the two following
cone decompositions:

\[\epsfbox{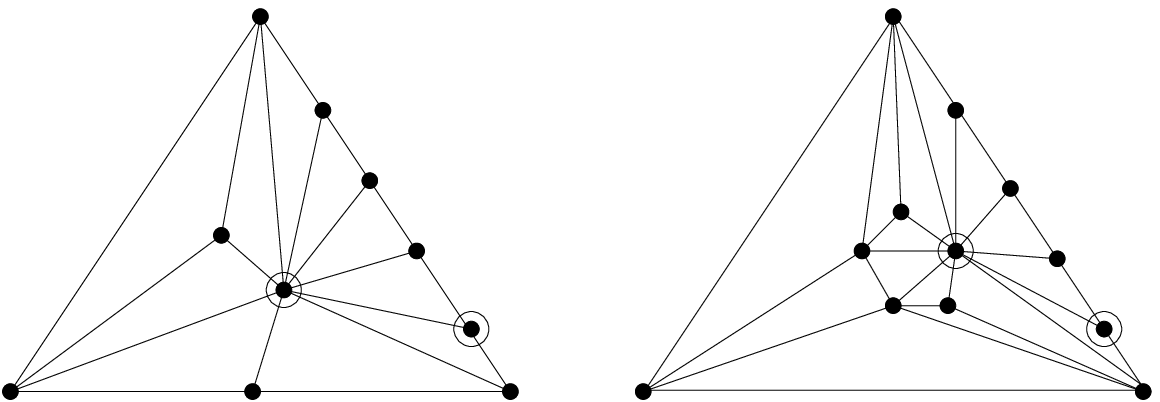}\hspace*{1cm}\]

From the singular curves we get $1+4=5$ components, from each $\fZ_{10}$
point an additional one and at each $\fZ_{15}$ point, we get four further
compoenents. Hence the singular fiber has a total of 15 components, in
addition to $\Theta$, the proper transform of $C_{\infty}$. 

\begin{figure}
\[ \unitlength1cm\begin{picture}(10,9.2)
\put(-3,0){\epsfbox{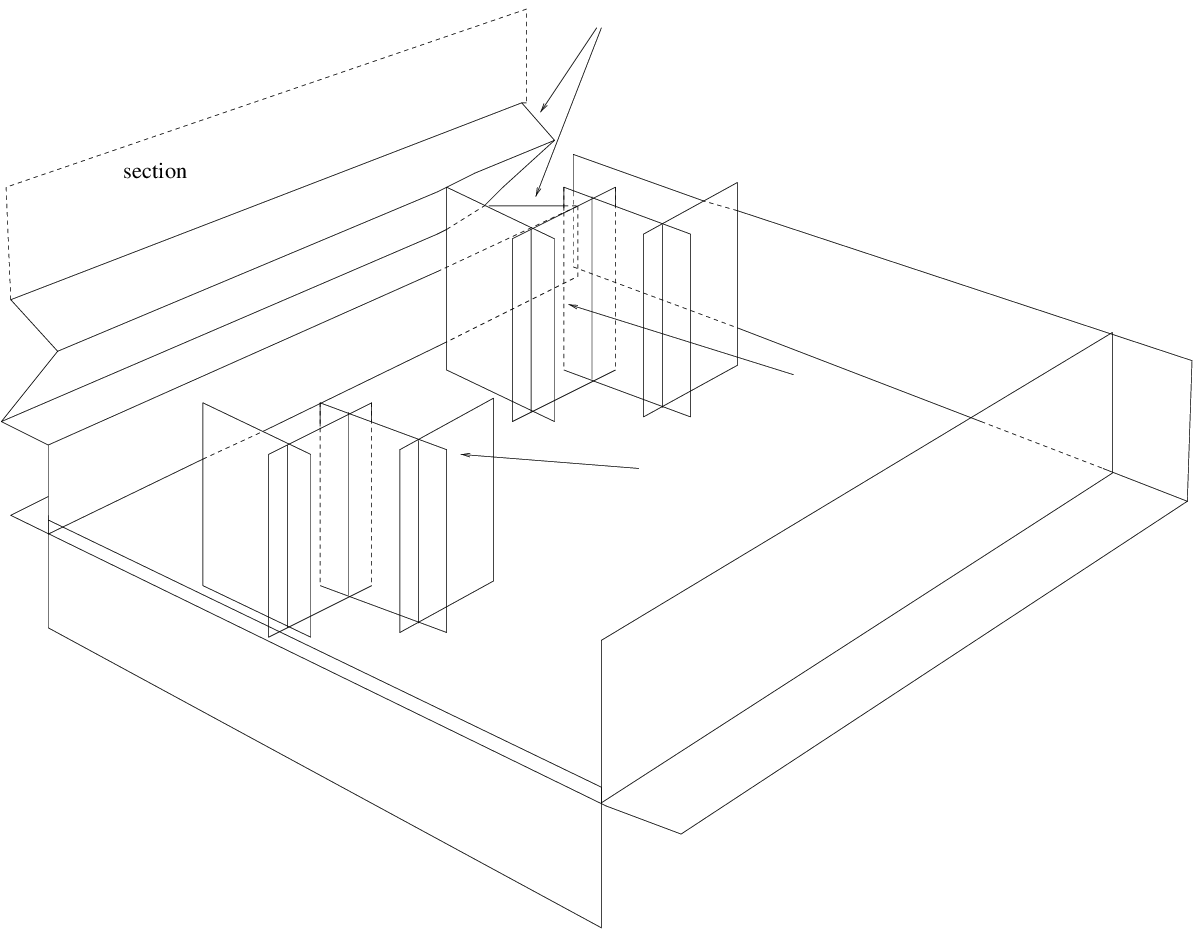}}
\put(3,9.2){$\Theta_{1,1}-\Theta_{1,4}$} 
\put(3.6,4.5){$\Theta_{3,1}-\Theta_{3,4}$}
\put(4.6,5.3){$\Theta_{4,1}-\Theta_{4,4}$}
\put(2.5,.5){$\Theta_1$}\put(7.5,5.2){$\Theta_2$}
\put(4,2.8){$\Theta_{2,1}$} \put(2.5,2){$\Theta$}

\put(4.5,8.5){The singular fiber is the union of 16 components}
\put(4.5,8){$\Theta+\Theta_{1,1}+\cdots+\Theta_{1,4}+\Theta_{2,1}
  +\Theta_1+\Theta_2+$}
\put(5,7.5){$\Theta_{3,1}+\cdots +\Theta_{3,4}+\Theta_{4,1}+\cdots + 
\Theta_{4,4}$} 
\put(6,7){where $\Theta$ is the proper transform of $C_{\infty}$, and}
\put(7.5,6.5){$\Theta_{1,i},\ i=1,\ldots, 4$ are the}
\put(9,6){components resolving the}
\put(9.5,5.5){$\fZ_5$ curve, $\Theta_{2,1}$ resolving}
\put(9.5,5){the $\fZ_2$ curve}
\put(9.5,4.5){$\Theta_{i,j}, i=3,4,\ j=1,\ldots, 4$}
\put(9.5,4){are the divisors }
\put(9,3.5){resolving the $\fZ_{15}$ points,}
\put(8.5,3){$\Theta_i,\ i=1,2$}
\put(8,2.5){resolving the $\fZ_{10}$ points.}

\end{picture}
\]
\caption{\label{figure7} --- The singular fiber of type IX$^{***}_1$ ---}
\end{figure}

After resolution of singularities, the singular fiber looks as in Figure
\ref{figure7}. 

\subsection{Analogs of Kodaira's type II$^*$ fiber}
{\it Case} IX$^{****}_1$:

The ambient space is $\fP_{(1,6,14,21,42)}$, which has the following
singular locus:
\begin{itemize}\item $\gS_1=\{z_1=z_2=0\}_{\fZ_7}$, \item
  $\gS_2=\{z_1=z_3=0\}_{\fZ_3}$, \item $\gS_3=\{z_1=z_4=0\}_{\fZ_2}$, 
  \item $\gS_4=(\gS_1\cap \gS_2)_{\fZ_{21}}$, \item $\gS_5=(\gS_1\cap
    \gS_3)_{\fZ_{14}}$,  \item $\gS_6=(\gS_2\cap \gS_3)_{\fZ_6}$. 
  \end{itemize}
The intersections $\gS_i\cap X$ are all rational curves, which meet two at
a time. We have the following configuration in $C_{\infty}$:
\[\unitlength1cm
\begin{picture}(10,3.5)
\put(0,0){ \epsfbox{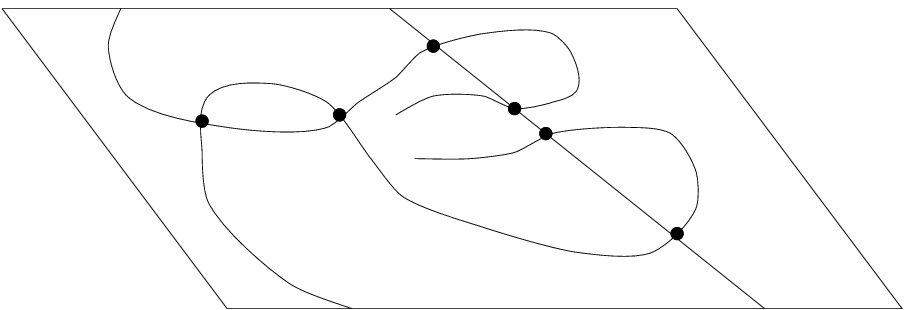} }
\put(1.8,2.3){$\fZ_6$}\put(3.3,2.3){$\fZ_6$}
\put(1.2,3.3){$\gS_2 -\ \fZ_3$} \put(3,-.3){$\gS_3 -\ \fZ_2$}
\put(7,-.3){$\gS_1 -\ \fZ_7$}
\put(7.5,1){$C_{\infty}$} \put(5.4,1.4){$\fZ_{14}$}
\put(6.1,.8){$\fZ_{14}$}
\put(3.7,2.8){$\fZ_{21}$} \put(5,2.3){$\fZ_{21}$}

\end{picture}
\]
\begin{figure}

\[ \unitlength1cm \begin{picture}(10,5) \put(-3,0){\epsfbox{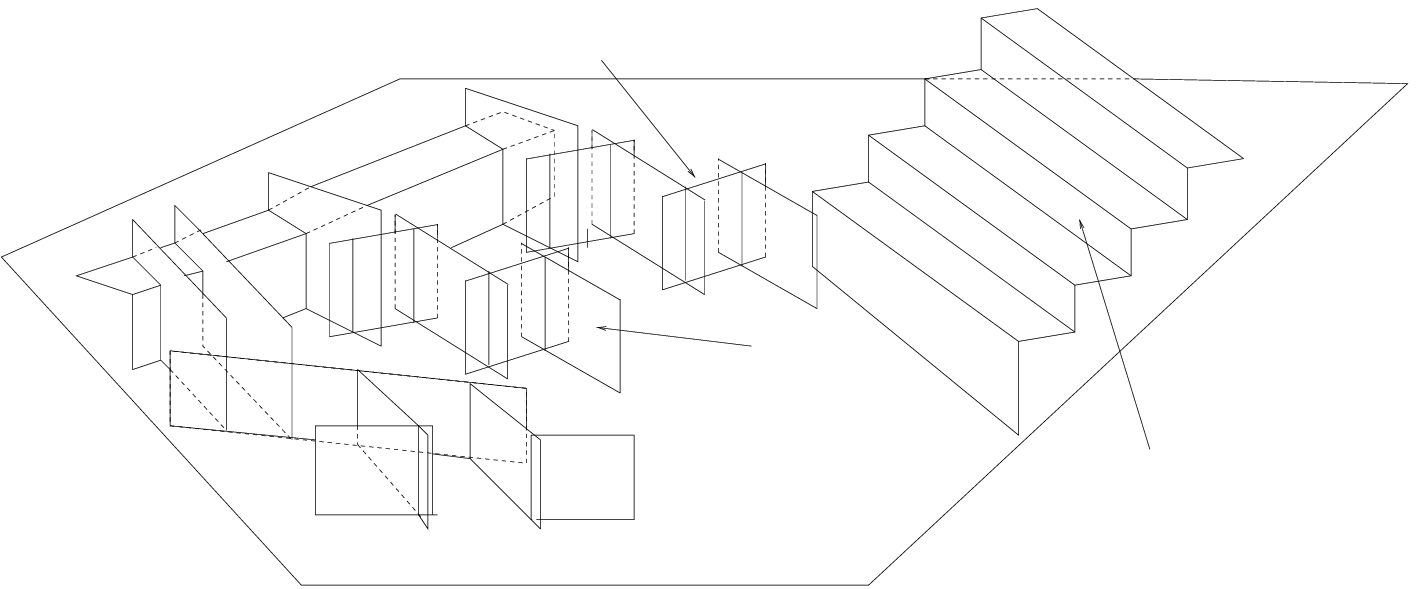}}
 \put(2.7,1){$\Theta_{5,2}$}\put(4.5,1){$\Theta$}
 \put(1.85,1.4){$\Theta_{5,1}$} \put(.5,1){$\Theta_{5,4}$}
 \put(.7,1.7){$\Theta_{5,3}$} \put(0,1.8){$\Theta_{3,1}$}
\put(-.7,2.2){$\Theta_{4,1}$} \put(-1.4,2.7){$\Theta_{4,2}$}
\put(1.6,3.9){$\Theta_{2,1}$} \put(.7,4.1){$\Theta_{2,2}$}
 \put(3.8,2.2){$\Theta_{4,1}-\Theta_{4,5}$}
\put(2.5,5.7){$\Theta_{4,6}-\Theta_{4,10}$}
\put(7.9,1){$\Theta_{1,1}-\Theta_{1,6}$}
\put(7.9,.5){$\Theta_{1,6}$ is a section, not a component}

\end{picture}
\]
\caption{ {\label{figure8}--- The singular fiber of type IX$^{****}_1$ ---} }

\medskip 
{The singular
 fiber is a sum of 26 components. $\Theta$ is the proper transform of
 $C_{\infty}$, and for the exceptional divisors we have choosen the
 notation so that the divisors $\Theta_{i,j}$ resolve the singular locus
 $\gS_j$. There are then $5+2+1+10+4+2 = 25$ components of the various loci. }
\end{figure}

\smallskip
The $\fZ_6$ points are resolved precisely as in the cases above; there is
one exceptional divisor, which this time is a component of the singular
fiber, as it is not contained in all fibers, but only in $C_{\infty}$. It
remains to resolve the $\fZ_{14}$ and $\fZ_{21}$ points. They are of types
${1\over 14}(1,6,7)$ and ${1\over 21}(1,6,14)$, and we have $d_1=1,\
d_2=2,\ d_3=7$ in the first case and $d_1=1,\ d_2=3,\ d_3=7$ in the second
case. By our formula above, this means we have $v=13$ and $v=17$,
respectively, leading to 3 and 6 inner vertices, respectively. We choose
the following cone decompositions:

\[ \epsfbox{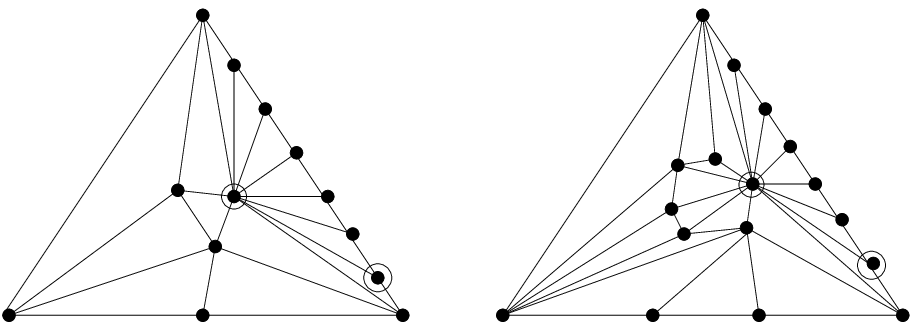} \]

Again the vertices corresponding to the sections of the fibration are
circled. Just as above, from this we can without difficulty derive the
singular fiber. It will look as in Figure \ref{figure8}.

\smallskip\noindent{\it Case} XII$^{***}_3$:

The ambient projective space is $\fP_{(1,6,42,98,147)}$ with singular locus

\begin{figure}
\[\epsfbox{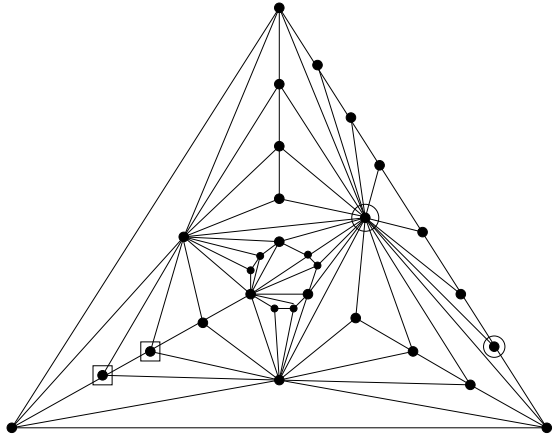} \]
\caption{\label{figurecone} The cone decomposition for the case
  XII$^{***}_3$} 
\end{figure}

\begin{itemize}\item $\gS_1=\{z_1=z_2=0\}_{\fZ_7}$, \item
  $\gS_2=\{z_1=z_4=0\}_{\fZ_3}$, \item $\gS_3=(\gS_1\cap \gS_2)_{\fZ_{21}}$,
  \item $\gS_4=\{z_1=z_2=z_3=0\}_{\fZ_{49}}$. 
  \end{itemize}
The intersections $\gS_i\cap X,\ i=1,2$ are both rational curves, which
meet in a single point. Furthermore, $\gS_4\cap X$ consists also of a
single point. One sees easily that the $\fZ_{21}$ point is resolved exactly
as above in the case IX$^{***}_1$. It remains to resolve the $\fZ_{49}$
point. Since this lies on a $\fZ_7$ curve, the singularity is of the type
${1\over 49}(1,6,42)$, and we have $d_1=d_2=1,\ d_3=7$. The number of
vertices is then ${51+9\over 2} = 30$, of which $3+6$ are on the
boundary. It follows that we have to include $21$ inside vertices. We take
the cone decomposition of Figure \ref{figurecone}.

\begin{figure}
\[
\unitlength1cm \begin{picture}(10,8.5)
\put(-2,0){\epsfbox{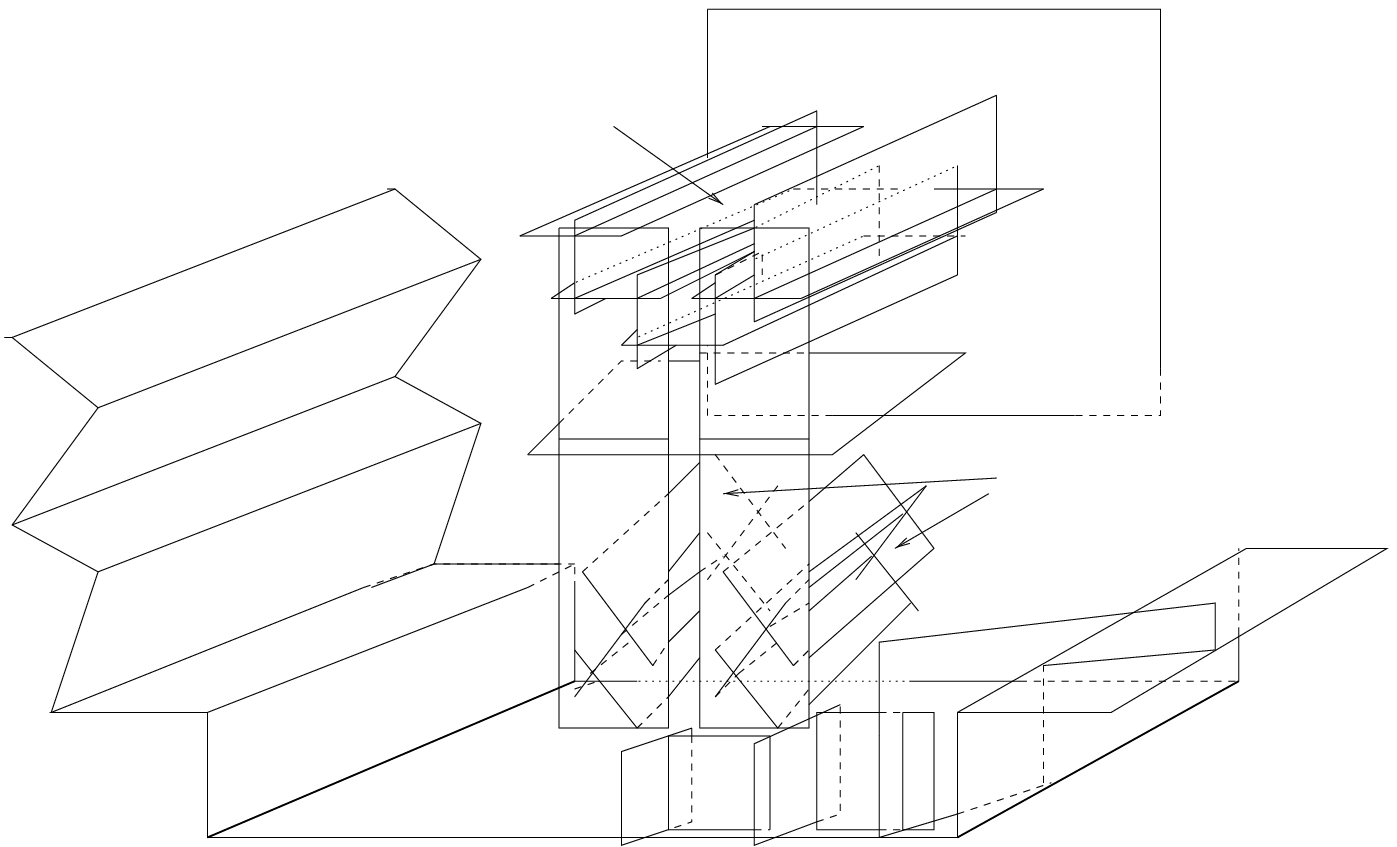}}
\put(2,.5){$\Theta$}
\put(1.5,1.3){$\Theta_{1,1}$}
\put(.8,2){$\Theta_{1,2}$}
\put(.7,2.8){$\Theta_{1,3}$}
\put(0,3.5){$\Theta_{1,4}$}
\put(0,4.4){$\Theta_{1,5}$}
\put(0,5.8){$\Theta_{1,6}$ (a section)}
\put(8.5,5){$\Theta_{4,9}$}
\put(3,7.5){$\Theta_{4,1}-\Theta_{4,8}$}
\put(8,3.8){$\Theta_{4,12}-\Theta_{4,14}$}
\put(8,3.3){and $\Theta_{4,15}-\Theta_{4,17}$}
\put(6.4,4.6){$\Theta_{4,18}$}

\put(4.35,.5){$\Theta_{3,1}$}
\put(5.1,.7){$\Theta_{3,2}$}
\put(5.8,.3){$\Theta_{3,3}$}
\put(6.4,.5){$\Theta_{3,4}$}
\put(7.4,1.6){$\Theta_{3,5}$}
\put(3.8,3.3){$\Theta_{4,10}$}
\put(5.3,3.3){$\Theta_{4,11}$}
\put(8.7,1.1){$\Theta_{2,1}$}
\put(10.3,2.5){$\Theta_{2,2}$}

\end{picture}
\]
\caption{\label{figure9}--- The singular fiber of type XII$_3^{***}$ ---}

\medskip
The singular fiber is the union of 33 components. $\Theta$ is the proper
transform of $C_{\infty}$, and again the components $\Theta_{i,j}$ are the
exceptional divisors resolving $\gS_i$. Of the 20 components resolving the
$\fZ_{49}$ point, we have not drawn two of them, which correspond to the
vertices labeled with a square in the above cone decomposition, as they
would have cluttered up the picture too much. 
\end{figure}

The singular fiber will look quite a bit like that of type XII$_3^{**}$,
but will have nine additional components. This will look as displayed in
Figure \ref{figure9}.

\subsection{The other cases}

We now consider the last three cases in Table 3. The difference between
these and the above cases is that there are now two bad fibers of $^*$
type. These are the fibers which are the total transforms of the surfaces
$C_{\infty}$ and $C_0$. The fiber $C_0$ is now also singular because of the
fact that the first weight is no longer unity. A detailed analysis is not
necessary. Consider case 13, i.e., the projective space
$\fP_{(2,3,5,10,10)}$. The surface $C_{\infty}$ is
$\fP_{(2,5,10,10)}[30]\cong\fP_{(1,1,1,1)}[3]$, a cubic surface, and the
two singular curves on it are elliptic and meet in three points. Without
difficulty we recognize the singular fiber of type IV$^{*}_1$. The surface
$C_0$ is $\fP_{(3,5,10,10)}[30]\cong \fP_{(1,2,2,3)}[6]$, and we recognize
the singular fiber of type IV$^{**}_1$. Note that the singular curve
$\gS_1-\{z_1=z_2=0\}$ is $\fZ_5$, and yields upon resolution four components,
one of which is a section of the fibration. The other three split; one of
the components belongs to the fiber at $\infty$, while the other belongs to
the fiber of type IV$^{**}_1$, and indeed, the first has one, the second
has two such components. Similarly, one can check that the exceptional
divisors which lie over the singular points split, one component is a
section, the others belong to one or the other fiber.  
The same methods apply to the remaining cases, and the results are listed
in Table 3.


\begin{thebibliography}{M}

\bibitem{aspin} P.\ Aspinwall, {\it (2,2)-superconformal theories near
    orbifold points}, Comm.\ Math.\ Phys. {\bf 128} (1990), 593-611.

\bibitem{Y} A.\ Dimca, {\it Singularities and coverings of weighted
    complete intersections}, J.\ reine und angew.\ Math.\ {\bf 366} (1986),
    184-193.

\bibitem{dolg} I. Dolgachev, {\it Weighted projective spaces}, in ``Group
  actions and vector fields'' (Vancouver, B.C., 1981), Springer Lecture
  Notes in Mathematics {\bf 966}, Springer,
  1982, pp. 34-71.

\bibitem{GRY} B.\ Greene, S.\ Roan \& S.-T.\ Yau, {\it Geometric
    singularities and spectra of Landau-Ginzburg Models}, Commun.\ Math.\
    Phys.\ {\bf 142} (1991), 245-259.

\bibitem{I} B.\ Hunt \& R.\ Schimmrigk, {\it On K3-fibered Calabi-Yau
    threefolds I, the twist map}, MPI preprint \# 48/1998, this Journal.
 
\bibitem{RY} S.\ Roan, S.-T.\ Yau, {\it On Ricci flat 3-fold}, Acta
  Mathematica Sinica, New Series {\bf 3} (1987), 256-288.
\end{thebibliography}
\end{document}